\documentclass[11pt,a4]{amsart}
\usepackage{parskip}
\usepackage{amsfonts,amsmath,amssymb,bbm}
\usepackage{verbatim}
 \usepackage{setspace}
\usepackage{color}
\usepackage[colorlinks=true]{hyperref}

\makeatletter
\renewcommand{\section}{\@startsection{section}{1}{0pt}{20pt}{6pt}{\large\bfseries}}
\makeatother

\usepackage{enumerate}
\usepackage{pdfsync}
\usepackage{color}
\usepackage{dsfont}

\numberwithin{equation}{section}

\theoremstyle{plain}
  \newtheorem{theorem}{Theorem}%[section]
  \newtheorem{proposition}[theorem]{Proposition}
  \newtheorem{lemma}[theorem]{Lemma}
   \newtheorem{corollary}[theorem]{Corollary}

\theoremstyle{definition}
\newtheorem{definition}[theorem]{Definition}

\theoremstyle{remark}
\newtheorem{remark}[theorem]{Remark}

\newcommand{\E}{\mathbb{E}}

\newcommand{\R}{\mathbb{R}}
\newcommand{\C}{\mathbb{C}}
\newcommand{\N}{\mathds N}
\newcommand{\Z}{\mathbb Z}
\newcommand{\D}{\mathds{D}}
\newcommand{\Pp}{\mathds{P}}
\newcommand{\Dt}{\mathds{D}_2}
\newcommand{\Dp}{\mathds{D}_P}
\newcommand{\De}{\mathds{D}_L}
\newcommand{\Nd}{\mathds{N}_{\D}}
\newcommand{\Be}{\mathds{B}}
\newcommand{\Bd}{\mathds{B}_{\D}}

\newcommand{\WPn}{W_{\Psi}(n+1)}

\title{A new class of solutions to the {van} Dantzig problem, the Lee-Yang property, and the Riemann hypothesis}

\author{T.~Konstantopoulos}
\address{Department of Mathematical Sciences, University of Liverpool, Liverpool, L69 7ZL, UK}
\email{takiskonst@gmail.com}

\author{P.~Patie}\thanks{The second author is grateful to Mateusz Kwa\'snicki for providing the reference \cite{Wid} which contains a closure property of P\'olya frequency densities.}
\address{School of Operations Research and Information Engineering, Cornell University, Ithaca, NY 14853.}
\email{pp396@cornell.edu}

\author{R.~Sarkar}
\address{School of Operations Research and Information Engineering, Cornell University, Ithaca, NY 14853.}
\email{rs2466@cornell.edu}

\usepackage[usenames,dvipsnames]{xcolor}
\newcommand{\red}{\color{red}}

\usepackage{etoolbox}
\makeatletter
\usepackage{cancel}
\newcommand\Ccancel[2][black]{
    \let\OldcancelColor\CancelColor
    \renewcommand\CancelColor{\color{#1}}
    \cancel{#2}
    \renewcommand\CancelColor{\OldcancelColor}
}
\renewcommand\CancelColor{\red}

\usepackage{bm}
\def\LPnew{\mathds{LP}_{\hspace{-0.06cm} e}} % even F in LP with F(0) =1

\begin{document}

\maketitle

\begin{abstract}
The purpose of this paper is to carry out an in-depth analysis of the intriguing  van Dantzig
problem which  consists on characterizing the set $\D$ of  analytic characteristic functions $\mathcal{F}$  which remains stable by the action of the mapping $V\mathcal{F}(t)=1/\mathcal{F}(it)$, $t\in\R$. % is also a characteristic function.
 We start by observing that the celebrated Lee-Yang property,
 appearing in statistical mechanics and quantum field theory,
and the Riemann hypothesis can be both rephrased in terms of the van Dantzig problem, and, more specifically, in terms of the set $\De \subset \D$  of real-valued characteristic  functions that belong to the Laguerre-P\'olya class. Motivated by these facts, we proceed by identifying several non-trivial closure properties of the set $\D$ and $\De$.  This not only revisits but also,  by means of probabilistic techniques,   deepens the fascinating studies of the set of even characteristic functions in the  Laguerre-P\'olya class  carried out by P\'olya \cite{Polya27}, de Bruijn \cite{DB},  Lukacs  \cite{Lukacs}, Newman \cite{New1} and more recently by Newman and Wu \cite{New}, among others.  We continue by providing  a new class of entire functions that belong to the set $\D$
but not necessarily to $\De$, offering the first examples  outside the set $\De$.
This class, which is derived from some entire functions introduced by the second author in \cite{Patie}, is in bijection with a subset of continuous negative-definite functions and includes several notable generalized hypergeometric type functions. Besides identifying the characteristic functions,
we also manage to characterize  the pair of the corresponding van Dantzig random variables revealing that  one of them is  infinitely divisible. Finally, we investigate the possibility that the Riemann $\xi$ function belongs to this class.

\end{abstract}


%%%%%%%%%%%%%%%%%%%%%%%%%%%%%%%%%%%%%%%%%%%%%%%%%%%%%%%%%%%%%%%%

%%%%%%%%%%%%%%%%%%%%%%%%%%%%%%%%%%%%%%%%%%%%%%%%%%%%%%%%%%%%%%%%
%$\vspace{2pt}$
%\\

%{\footnotesize\tt
%\tableofcontents
%}

%%%%%% INSERTED \input{Intro-Properties.tex}  HERE  %%%%%%%%%

\section{Introduction}\label{introduction}%to the van Dantzig problem
%\subsection{van Danzig pairs}
In his seminal paper \cite{Lukacs}, E.\ Lukacs provided  comprehensive
and fascinating insights towards the characterization of the class of functions
\begin{equation}    \label{eq:defD}
\D= \left\{\mathcal{F} \in \Pp_+; \:  V\mathcal{F} \in \Pp_+    \right\}
\end{equation}
where
\[
V\mathcal{F}(t)=\frac{1}{\mathcal{F}(it)},t \in \R,
\]
and $\Pp_+$ stands for the set of continuous bounded positive-definite functions,
normalized to be $1$ at the origin. Note that by Bochner's theorem,
 $\Pp_+$ is
the set of characteristic functions of real-valued  random variables, see e.g.~\cite{SchillingSongVondracek10}.
 This question  was posed,
as a prize-winning problem, by David van Dantzig, a Dutch algebraic topologist,
in Nieuw Archief voor Wiskunde.

If $\mathcal F \in \mathds D$ then the function $V\mathcal{F}$
is called the (van Dantzig) reciprocal of $\mathcal F$, and also belongs to $\D$.
 This last fact  entails, see Proposition \ref{prop:propD}, that
any $\mathcal{F} \in \D$
admits an analytic extension on a cross section of the complex plane
including the imaginary line in its interior.
We also let
\begin{equation}\label{eq:defD}
  \Dt= \left\{[\mathcal{F},V\mathcal F] \in \Pp_+\times \Pp_+ ; \: \mathcal{F}(it) \cdot V\mathcal F(t)=1 \textrm{ for all } t \in \R\right\}.
\end{equation}
 We say that $\mathcal F$ is a van Dantzig function
and that $[\mathcal{F},V\mathcal F]$ is a van Dantzig pair. Note that the mapping $V: \mathds D \to \mathds D$ is a multiplicative involution, that is,
$V \mathcal F_1 \mathcal F_2=V \mathcal F_1 V \mathcal F_2$ and $V \circ V\mathcal F =\mathcal F$, see Proposition \ref{prop:invol} below.

The first  historical instances of non-trivial van Dantzig pairs are
\begin{equation}\label{eq:ex}
  \left[ \cos t, \frac{1}{\cosh t}\right], \left[ \frac{\sin t}{t}, \frac{t}{\sinh t}\right] \textrm{ and }\left[ e^{-\frac{t^2}{2}}, e^{-\frac{t^2}{2}}\right].
\end{equation}
We notice that, in the last example,
\[
\mathcal{F}_{\textrm{N}}(t)=e^{-\frac{t^2}{2}}
\]
is the characteristic function of a standard normal random variable for which
it is immediate
that $V\mathcal F_{\textrm{N}}(t)=1/\mathcal{F}_{\textrm{N}}(it)=\mathcal{F}_{\textrm{N}}(t)$ and hence $\mathcal F_{\textrm{N}} \in \D_S$, where
\[
\D_S = \{\mathcal F \in \Pp_+;\: \mathcal F= V\mathcal F\}
\]
is the invariant set  of $V$ or self-reciprocal elements of $\D$. It is easy to see that, for any $\mathcal F \in \D$, the mapping
${t\mapsto\mathcal{F}(t)} V\mathcal F(t)$
is in $\D_S$, meaning that the  set $\D_S$ contains many more elements than $\mathcal F_{\textrm{N}}$.
However, the following fact, due to Lukacs, offers an original characterization of the characteristic function $\mathcal F_{\textrm{N}}$
%\begin{enumerate}[(i)]
%% \item If  $\mathcal{F} \in  \De$ and $t\mapsto  \log \mathcal{F}(it) \in \N(\R)$.
%\item \label{it:fn}
%$\mathcal{F}_N \in \D_S$ and $\{ [\mathcal{F}_1,\mathcal{F}_2] \in \D_2; \:  -\log %\mathcal{F}_j \in \N(\R), j=1,2\}=\{[\mathcal{F}_N,\mathcal{F}_N]\}$.
%\end{enumerate}
\[
\big\{ [\mathcal{F}_1,\mathcal{F}_2] \in \D_2; \:  -\log \mathcal{F}_j \in \N(\R),\, j=1,2\big\}
=\big\{[\mathcal{F}_{\textrm{N}},\mathcal{F}_{\textrm{N}}] \big\}
\]
where $\mathds N(\R)$ is the set of negative definite functions on $\R$,
 see  \cite[Chapter 4]{SchillingSongVondracek10} for more information on this set. In other words, if  $\mathcal F$ and $V \mathcal F$ are
  infinitely divisible then they are identical and both equal to
$\mathcal{F}_{\textrm{N}}$.

 Regarding the first example in \eqref{eq:ex}, it is immediate that the mapping $t\mapsto \cos t$ is the characteristic function of
a random variable taking values $\pm1$ with equal probability, whereas
$t\mapsto(\cosh t)^{-1}$ is the characteristic function of a random variable with density
$(2\cosh(\pi x/2))^{-1}$,  $x\in \R$.  For the second example, the mapping $t\mapsto \sin t/t$ corresponds to the characteristic function of a uniformly distributed
random variable on the interval $[-1,1]$
and $t\mapsto t(\sinh t)^{-1}$   to the one of an absolutely continuous probability measure whose  density is $\frac{\pi}{4}(1-\tanh(\pi x/2))^2,  x\in \R$.

Another classical example, discussed at the end of Lukacs's paper, and
presented in a more general form  in \cite{Jurek,Hinds},  is expressed in terms of the entire functions
\begin{equation}\label{eq:bes}
\mathcal{J}_{\nu}(t)=\Gamma(\nu+1)t^{-\nu} J_{\nu}\left( t\right), \quad
\mathcal{I}_{\nu}(t)=\mathcal{J}_{\nu}(it),
\end{equation}
where $J_{\nu}$ is the Bessel function of the first kind of index $\nu$, and given by  the pair
\begin{equation}\label{eq:ex2}
 \left[ \mathcal{J}_{\nu},\, \frac{1}{ \mathcal{I}_{\nu}}\right] \in \Dt,
\quad \nu  > - \frac12.
\end{equation}
The corresponding pair of random variables are described in details in \cite{Jurek}, see also \eqref{eq:defJ} and Lemma \ref{lem:BT}. This set of  examples turns out to be the only canonical solutions to the van Dantzig problem available in the literature. Note that they  all belong to the subclass $\D_L \subset \D$, originally identified by  Lukacs  \cite{Lukacs}, which stands for
the set of even entire characteristic functions in the Laguerre-P\'olya class.
That is, entire functions which are locally the limit of a series of polynomials whose roots are all real, see Section \ref{sec:DL} for further discussion on this class.
Throughout, we shall provide several ways of generating new instances in $\D_L$ and also present  new subclasses of $\D$.

In this paper, we start by identifying a connection between two fundamental problems in mathematics and the subclass $\mathds D_L \subseteq \D$.
First, we explain how the Riemann hypothesis can be equivalently formulated by the membership of
the Riemann $\xi$ function to the class $\mathds D_L$, something which, in a different context,  was observed by Roynette and Yor \cite{Roy}. Similarly, the celebrated Lee-Yang property in statistical physics, discovered first
by Lee and Yang  in connection to the Ising model on a finite lattice,
is also equivalent to  the requirement that the partition function, viewed as a function on the imaginary line, belongs to $\mathds D_L$, see \cite{LYs} for a thorough account on this topic. This connection relies on the theory of P\'olya frequency functions developed by  Schoenberg \cite{Schoenberg}, see Section \ref{sec:DL}.

In view of its importance, we first aim  to develop an in-depth analysis of the subclass $\D_L$. On the one hand, we adapt  several substantial results that one can find in the  number theory and statistical mechanics literature  to provide new information about the class $\D_L$. Conversely, by means of probabilistic techniques combined with the theory of entire functions, we also provide original closure properties of the class $\D_L$ which give new insights to the two aforementioned problems. For instance, we  find necessary conditions for the  product of two independent random variables to belong or to remain in $\D_L$, see Theorem \ref{thm:DM}. In the same vein, in Theorem \ref{thm:conv}, we   revisit, improve and extend to the class $\D$, a recent result due to Newman and Wu \cite{New} regarding  the closure property of this class under locally uniform convergence.

Another objective of this paper is to identify a new subclass, denoted by $\D_P$, of analytic characteristic functions that belong to the class $\mathds D$, that is, are solutions to the van Dantzig problem.
More specifically, to each Laplace exponent
$\Psi$ of a (possibly killed) spectrally negative L\'evy process which  is {{non-negative}}  on the
interval $[1/2, \infty)$, we associate  the function  $\mathcal J_\Psi \in \mathds D_P$ which is defined by
\[
\mathcal J_\Psi(t) = 1+\sum_{n=1}^\infty
\frac{(-1)^n}{\Psi(1) \cdots \Psi(n)} t^{2n}.
\]
Since $\lim_{u \to \infty}\Psi(u)= \infty$, we deduce easily  that
$\mathcal J_\Psi$ defines an entire function. We refer to sections \ref{intro part 2} and \ref{sec:PropJ} for  more information about the objects introduced here and this class of functions, whose analysis is intimately related to
the Wiener-Hopf factorization of the Laplace exponent
$\Psi$. We also point out that the entire function $\mathcal I_\Psi(t)=\mathcal J_\Psi(i\sqrt{t})$ was introduced by the second author in  \cite{Patie}, and has appeared in various mathematical contexts recently, in  probability theory  \cite{Patie1}, in the spectral theory of some non-self-adjoint operators \cite{Patie-Savov, Patie-Sarkar}, in the study of special functions \cite{Zhang} and a chapter is devoted to this class of function in the monograph \cite{Kyprianou06}. Note that when $\Psi(n)=n(\nu+n)$ is the Laplace exponent of a (scaled) Brownian motion with drift {{$\nu\geq -\frac12$}}, then $\mathcal{J}_{\Psi}$ boils down to the Bessel-Clifford  function
$\mathcal{J}_{\nu}$ as defined in \eqref{eq:bes}. The class $\D_P$  includes a wide range of other well-known special
functions such as several hypergeometric functions, the Mittag-Leffler functions, and
the Wright functions, among others, see Section \ref{sec:exa}.

Theorem \ref{thm:main} states that $\mathcal J_\Psi \in\mathds D$.  Its proof relies on identifying the two random variables
whose characteristic functions are $\mathcal J_\Psi$ and $1/\mathcal I_\Psi, \mathcal I_\Psi(t)=\mathcal J_\Psi(it)$, see Lemmas  \ref{lem:IJ} and \ref{lem:BT}, respectively. This is achieved by introducing a Markov operator, see \eqref{eq:MI}, that serves on the one hand to show that $\mathcal J_\Psi$ is the  characteristic function of the  product of two independent  random variables (one having  as characteristic function the Bessel function $\mathcal J_0$). On the other hand, it also turns out to be  an intertwining
operator between two Markov semigroups, see e.g.~\cite{Patie-Savov} for a review of this concept. From this fact, we  deduce first that $t\mapsto \mathcal I_\Psi(\sqrt{t})$ is an invariant function for one of these semigroups and then  its reciprocal is the Laplace transform of a positive random variable associated to the Markov process, from where we conclude by invoking an argument involving the Bochner subordination of a Brownian motion. Moreover, since $\mathcal{J}_{\nu} \in \D_L \cap \D_P$, it is natural to wonder whether $\D_P \subset \D_L$. This is a very delicate question as it is difficult, in general, to identify the location of zeros of a power series.   However, we manage to provide necessary conditions on $\Psi$ for $\mathcal{J}_{\Psi} \in \D_L$, that is the entire function has only real zeros. We also identify instances of entire functions in $\D_P$ which have non real zeros, revealing that $\D_L \subsetneq \D$, see Theorem \ref{thm:pick}.

The remaining part of the paper is organized as follows. In Section \ref{sec:DL}, we introduce the  Lukacs class $\D_L$ and present its connection with the Riemann hypothesis and the Lee-Yang property. Section \ref{sec:Luk} is devoted to a thorough analysis of the sets $\D$ and $\De$,
including some closure properties of $\De$ under various mappings,
and to the identification of  original ways to generate new elements in $\De$.
In Section \ref{intro part 2} we introduce the new subclass $\D_P$ and provide some interesting properties.
Finally, in Section \ref{sec:proof}, we collect the proofs of the results presented in the two previous sections.

\section{The Lukacs class $\De$}  \label{sec:DL}
%partition function, the Lee-Yang property and the Riemann hypothesis}
We start by recalling that an entire function is in the Laguerre-P\'olya class $\mathds{LP}$ if it is the
local uniform limit of a sequence of polynomials with real coefficients and real
zeros only. In fact, based on ideas of Laguerre, see \cite[pp.\ 167--177]{Oeuvres} and
\cite[Ch.\ 2]{Borel}, P\'olya and Schur \cite{PolyaSchur} showed that
$\varphi$ is in $\mathds{LP}$ if and only if
\begin{equation}
\label{LPclass}
\varphi(z) = K z^m e^{-c^2 z^2 + az} \prod_{k=1}^\infty
\left(1-\frac{z}{z_k}\right) \, e^{z/z_k}, \: z\in \C,
\end{equation}
for some $m \in \Z_+$, $K, c, a \in \R$, $z_k \in \R\setminus\{0\}$ such that
$\sum 1/z_k^2<\infty$, and, where $z_k$, $k \in \N$, are the nonzero zeros of the entire function $F$, arranged in order of nondecreasing modulus.

Lukacs was interested in even characteristic functions in the class $\mathds{LP}$
and made the observation that such functions are automatically in $\D$, something
that we explain in the sequel.

Let now
\begin{equation}
\LPnew=\{\varphi \in \mathds{LP}; \:\varphi \textrm{ is even and } \varphi(0)=1\}.
\end{equation}
 The representation \eqref{LPclass} immediately gives that
$\varphi \in \LPnew$ must be of the form
\begin{equation}
\label{eq:defPL}
%F \in \mathds{L} \iff
\varphi(z) = e^{-c^2 z^2} \prod_{k=1}^\infty
\left(1-\frac{z^2}{z_k^2}\right), \quad
c \in \R, \quad
z_k >0,\, \quad
\sum_k 1/z_k^2 < \infty.
\end{equation}
Indeed, the even functions in $\mathds{LP}$ are the functions of the
form \eqref{LPclass} with $m=0$, $K=1$ and zeros that are symmetrically
placed around $0$ on the real axis, resulting in precisely the form
\eqref{eq:defPL}. Next,  we mention that the order $\rho$ and the exponent of convergence $\varrho$ of an entire function $\varphi$ in the Laguerre-P\'olya  class are such that \[ 0\leq \varrho \leq \rho \le 2\] where we recall that
$\rho =  \overline{\lim}_{r \to \infty} \frac{\log\log \max_{|z|=r} |\varphi(z)|}{\log r}$
and
$\varrho= \inf\left\{\alpha>0;\, \sum_{k\geq1} |z_k|^{-\alpha}< \infty
\right\}$.
Here and throughout, we refer to
the monograph of Levin \cite{Levin96} for information related to entire functions.

The functions in the examples \eqref{eq:ex} and \eqref{eq:ex2} are  all
of the form \eqref{eq:defPL}.
%$\mathds{LP}_2$.
However, we emphasize that  if an entire function $\varphi$ is of the form
\eqref{eq:defPL} then it needs
not be the case that the mapping $t \mapsto \varphi(t)$  is positive-definite on $\R$.
Nevertheless, Schoenberg, following Hadamard,   proved that the reciprocal of a normalized Laguerre-P\'olya entire function (non-necessarily even) is in $\Pp_+$. In particular, we have the following.
\begin{theorem}	\cite[Theorem 1]{Schoenberg} \label{thm:sch}
If $\varphi$ is of the form \eqref{eq:defPL} then
\begin{equation}
\label{Zeno}
t\mapsto \frac{1}{\varphi(it)}
=e^{-ct^2} \prod_{k=1}^{\infty}\left(1+\frac{t^2}{z^2_k}\right)^{-1}
\end{equation}
is the characteristic function of a symmetric P\'olya frequency density $f_P$,
namely,
\[
\mathcal F_{f_P}(t) = \frac{1}{\varphi(it)} = \int_\R e^{itx} f_P(x) dx
\]
where $f_P$ is a symmetric probability density function on $\R$
such that, for all $n\in \N$,
$x_1< \cdots < x_n$, $y_1 < \cdots < y_n$, the determinant of
the matrix $[f_P(x_j-y_k)]_{j,k=1}^n$ is non-negative.
\end{theorem}
Moreover, the probability measure with density $f_P$ is infinitely
divisible, that is, equivalently, $t\mapsto -\log \mathcal F_P(it) \in \N(\R)$,
see Kwa\'snicki \cite[Proposition 5.3]{Mateu}.

We point out that an easy way to see why
\eqref{Zeno} is a characteristic function is by a probabilistic argument. Indeed, let $\rm{N}, Z_1, Z_2, \ldots$ be independent random variables where
$\rm{N}$ is a standard normal and each $\rm{Z}_j$ a standard Laplace random variable,
that is, it has density $\frac{1}{2} e^{-|x|}$, $x \in \R$.
Then, one observes that the random variable
\begin{equation}
\label{NZ}
{\rm{X}} = c \sqrt{2} \textrm{N} + \sum_{j=1}^\infty \frac{\textrm{Z}_j}{z_j}
\end{equation}
has characteristic function \eqref{Zeno}. Then, relying on Theorem \ref{thm:sch}, Lukacs observed the very useful fact that if an even characteristic function $\mathcal F$ is in the Laguerre-P\'olya class then it is necessarily solution to the van Dantzig problem. This leads us to introduce the following class.

\begin{definition} \label{def:DL}
%Equivalently, $\LPnew$ contains precisely functions of the form \eqref{eq:defPL}.
Let \[ \D_L=\LPnew \cap \Pp_+ \subset \D\]
be the Lukacs class of solutions to the van Dantzig problem.
\end{definition}

The examples in \eqref{eq:ex} and \eqref{eq:ex2} are not just in $\D$
but also in $\D_L$.
Further elements of $\D_L$ can be generated by means of the mappings described in Theorem \ref{it:dif} below, and, also from the  subclass of $\D$ that we introduce and study in Section \ref{intro part 2}, see Theorem \ref{thm:pick}.
We already mention that in Section \ref{sec:Luk} we provide some (partial) characterizations of the class $\D_L$.

The first appearance of this  class of Laguerre-P\'olya  characteristic functions
traces back to P\'olya \cite{PolyaSchur} who, motivated
by Riemann hypothesis, was interested in characterizing all (complex valued)  functions
$f$ on $\R$ such that the analytic extension of
\begin{equation}\label{eq:pol}
t\mapsto  \int_{\R}e^{itx}f(x)dx \textrm{ is an entire function with only real zeros}.
\end{equation}
We refer to the recent paper by Newman and Wu \cite{NewRH}
for an excellent account on  P\'olya's  approach
and to de Bruijn's fascinating contributions to this problem.
Theorems \ref{thm:p1}, \ref{thm:p2} and \ref{thm:p3}  below form
the adaptation of these results in the context  of the van Dantzig problem.
From that time onwards, the Laguerre-P\'olya class of characteristic functions
has become  ubiquitous and plays a central role in  various fields of mathematics. In what follows, we describe the connection between the class $\D_L$ and both the Riemann hypothesis, and, the Lee-Yang property that appears in
 some statistical mechanics models and in Euclidean quantum field theory.

\subsection{The Lukacs class and the Riemann hypothesis}
There is a fascinating literature describing the role played by the Riemann $\zeta$ function in probability theory, see the excellent papers \cite{Biane,BianeP}  and the references therein. In the spirit of the work of Roynette and Yor \cite[Th\'eor\`eme V.3.2]{Roy}, we now explain how the Riemann hypothesis can be formulated in terms of the van Dantzig problem. % and  %The main interest in this question, from analytic number theory, is its relation to
%the Riemann hypothesis.
Consider the Riemann $\zeta$ function
\[ \zeta(z) = \sum_{n=1}^\infty n^{-z},\: \Re z > 1,\] which, when extended meromorphically to the whole complex plane,
has a single simple pole at $1$ with residue $1$ and the following zeros:
the trivial ones located  at the negative even integers and
the nontrivial ones lying in the critical strip $0< \Re z < 1$.
The Riemann hypothesis states that all nontrivial zeros are located
on the critical line $\Re z = 1/2$.
It is well-known that the function
\[
\eta(z) = \frac{1}{2} z(z-1) \pi^{-z/2}  \Gamma(z/2) \zeta(z),\quad z \in \C,
\]
%\tkside{I rebaptised $\xi$ as $\eta$ and  $\Xi$ as $\xi$.  \mbreak
%I don't mind any notation. \mbreak It can be changed again.}
is an entire function (because $1-z$ cancels the pole of $\zeta$,
whereas the trivial zeros of $\zeta$ cancel the poles of the gamma function $\Gamma(z/2)$, located at the
same places). Hence the zeros of $\eta$ are the nontrivial zeros of $\zeta$.
Moreover,  it satisfies
\[
\eta(z) = \eta(1-z), \quad z \in \C.
\]
Performing an affine transformation on $\C$ so that the critical line maps onto the
real line, we  get the Landau function
\[
\xi(z)  = \eta(\tfrac12 + iz), \quad z \in \C.
\]
The functional equation above then becomes
%\tkside{Typo corrected.}
\[
\xi(z) = \xi(-z), \quad z \in \C.
\]
Using the standard integral expression for the gamma function and the
definition of the $\zeta$ function one can  show
that
\begin{equation}
\label{xitrans}
\xi(t) = \int_{-\infty}^\infty e^{itx} \Phi(x) dx, \: t\in \R,
\end{equation}
where
\[
\Phi(x)
= \sum_{n=1}^\infty \left(4\pi^2 n^4 e^{9x/2}-6\pi n^2 e^{5x/2}\right)
\,
e^{-\pi n^2 e^{2x}}.
\]
It is plain that $\Phi(x)>0$ for all $x \geq0$.
Since $\xi$ is even, one gets  that $\Phi(x)>0$
for all $x \in \R$.
Also, $\Phi$ is integrable as
 $\int_{-\infty}^\infty \Phi(x) dx = \xi(0)<\infty$,
%and, being integrable, we get, since $\xi(0)=\frac{1}{2}$,  that the function $2\Phi$
%\footnote{ Remark that
%$\xi(0) = \eta(1/2) = \Gamma(1/4)|\zeta(1/2)|/(8\pi^{1/4}) \approx 0.497$
%Delete this...}
we have that the function $\Phi/\xi(0)$
is the density  of a symmetric real-valued random variable
and $\xi/\xi(0)$ is its characteristic function.
We have the following.
\begin{theorem}
\label{thm:RH}
The function $\xi/\xi(0)$ is in $\De$ if and only if the Riemann hypothesis holds.
\end{theorem}
\begin{proof}
If the function $\xi/\xi(0)$ is in $\De$ then it has only real zeros which means that
$\eta$ has all its zeros on the critical line. Hence the Riemann hypothesis holds.
If the Riemann hypothesis holds then $\xi/\xi(0)$ is an even entire characteristic
function with real zeros only. Using the result of  Proposition  \ref{DLD} below,
we have that $\xi/\xi(0)$ is in $\De$.
%\tkside{Please check the converse part.}
\end{proof}
The previous result combined with  Theorem \ref{thm:sch} yields this reformulation.
\begin{corollary}
The function $t \mapsto \xi(0)/\xi(it)$ is the characteristic function of a symmetric P\'olya frequency function  if and only if the Riemann hypothesis holds.
\end{corollary}
Further connections between the Riemann $\xi$ and  the van Dantzig problem will be discussed in Section \ref{sec:vDR}.
%Since $\eta(\frac12+it)/\eta(1/2)$ is a characteristic function,
%the Riemann hypothesis is implied by the statement that this function
%is in the Laguerre-P\'olya class.

\subsection{The Lukacs class and the Lee-Yang property} \label{sec:LY}
Entire characteristic functions with only real zeros  appear naturally in various
models of statistical mechanics and quantum field theory.
We briefly state the examples taken from the excellent paper  \cite{NewRH}. In their pioneering works,
Lee and Yang \cite{LY52,YL52} considered the Ising model in the presence of an external magnetic field and discovered that the zeros, in the magnetic field variable,
of the partition function are purely imaginary. This is equivalent to the following.
Let $\bm \mu$ denote the  probability measure on $\{-1,1\}^N, N\in \N,$ defined by
\[
\bm \mu(\bm x) = K e^{\sum_{j,k=1}^N J_{j,k} x_j x_k}, \quad \bm x=(x_1,\ldots,x_N),
\]
where, here and below, $J_{j,k} \ge 0$ and $J_{j,k} = J_{k,j}$ for all $j,k$ and $K$ is a positive normalizing constant. Then the partition function
\[
 \mathcal P_{\bm \lambda,\mu}(z) = \sum_{\bm x \in \{-1,1\}^N}
e^{ z \bm \lambda \cdot \bm x} \bm \mu(\bm x), \quad
\bm \lambda=(\lambda_1, \ldots \lambda_N) \in \R_+^N,
\]
is an entire function whose zeros   are all purely imaginary.
This was discovered by Lee and Yang and refer to it as \emph{the Lee-Yang property}.
We emphasize that the location and distribution of the zeros of the partition function are useful to determine substantial properties
of the underlying physical system such as phase transitions,  the infinite volume limit and existence of a mass
gap under an external magnetic field.
In relation to the van Dantzig problem, one observes that the mapping $ t \mapsto \mathcal P_{\bm \lambda,\bm\mu}(it)$
is the characteristic function of the random variable
\[
\bm \lambda \cdot \bm x
= \sum_{j=1}^N \lambda_j x_j
\]
under the measure $\bm \mu$. By Theorem \ref{thm:sch}, the Lee-Yang property is then equivalent to
the statement that this characteristic function is in $\mathds D_L$.
Generalizing the Lee-Yang result, Simon and Griffiths \cite{SG73} showed that if $\mu_0$
is a symmetric probability measure\footnote{In fact, a signed measure is
also allowed in   \cite{SG73}.}
on $\R$ such that
\[ \int_\R e^{b x^2} \mu_0(dx) < \infty \textrm{ for all } b \in \R \textrm{ and } \mathcal F_{ \mu_0}(z)\neq 0  \textrm{ for all } \Im(z)<0 \]
then, to the probability measure
\[
\bm \mu_\beta(d\bm x) = K  e^{\beta \sum_{j,k} J_{j,k} x_j x_k} \prod_{k=1}^{N}\mu_0(dx_k)
 \quad
\bm x =(x_1, \ldots, x_N) \in \R^N,
\]
there corresponds the partition function
\[
\mathcal P_{\bm \lambda,\bm \mu_\beta}(z) =\int_{\R^N}
e^{z\bm \lambda \cdot \bm x} \bm\mu_\beta(d\bm x), \quad \bm \lambda=(\lambda_1, \ldots,
\lambda_N) \in \R_+^N,
\]
with the property that, for all $\beta\ge 0$,
the zeros of the entire function $z \mapsto \mathcal P_{\bm \lambda,\bm \mu_\beta}(z)$ are all purely imaginary.
Hence, the mapping $t\mapsto \mathcal{P}_{\bm \lambda,\bm \mu_\beta}(it)$  is in $\mathds D_L$.
We also point out that these Lee-Yang type theorems arise in quantum field theory, and refer again to \cite{NewRH} for a more detailed description. We simply point out that an important  measure in this context is the measure
$\mu_0(dx) = e^{-a x^4 - b x^2} dx$, where $a>0$ and $x,b \in \R$. From Theorem \ref{thm:sch}, we deduce the following.

\begin{theorem} \label{thm:LY}
A partition function $\mathcal{P}$ has the Lee-Yang property if and only if $t\mapsto \mathcal{P}(it) \in \D_L$.
\end{theorem}
%\textcolor{blue}{Put a statement here saying that }

\section{Properties of the van Dantzig and  Lukacs classes}\label{sec:Luk}
%By items \eqref{it:even} and \eqref{it:comm} above, if $\mathcal G$ is the
%reciprocal of $\mathcal F$ then $\mathcal F$ is the reciprocal of $\mathcal G$.
We collect here several properties of both the classes $\D$ and  $\De$.  We start by presenting some basic properties that both classes share. Some of them were stated in Lukacs' paper \cite{Lukacs} without proof; for sake of conciseness, we indicate the main lines of proof.  Then, we move to the properties to the Lukacs class $\De$ which are of different types. Some results can be found  in the number theory or statistical mechanics literature in some form that we adapt, revisit or extend to identify new properties for the set $\De$. Moreover, based on ideas coming from probability theory,   we  also present  original and substantial results about this set, see Theorem \ref{thm:DM}, and, we improve, in Theorem \ref{thm:conv}, a very interesting closure property due to Newman and Wu. Let us start  with  the following simple but useful result which is a reformulation  of \cite[Theorems A.2.1 and A.2.2, page 335]{Linnik77}.
\begin{proposition} \label{prop:propD}
If $\mathcal{F} \in \D$ then $\mathcal{F}$ is a real and even function, meaning that the associated random variables are symmetric. Moreover, $\mathcal{F}$ admits an analytic extension to some cross $\{z\in \C; \:|\Im(z)|<z_{1}\}, \{z\in \C;  |\Re(z)|<z_{2}\}$ where $z_{1}>0,z_2>0$.  The same claims hold for $V\mathcal{F}$.
\end{proposition}

We also point out that since $\mathcal{F}$ is even, a theorem of Schoenberg %\cite{Sch1938}
 on positive-definite radial functions entails that the mapping
$t\mapsto \mathcal{F}(\sqrt{t})$ is completely monotone on $\R^+$, that is,
it is the Laplace transform of a non-negative Radon measure on $[0,\infty)$,
see e.g.~\cite{chafait}.%for a probabilistic proof,

We proceed with the following  the simple observation which follows readily from the previous Proposition since $V \circ V \mathcal{F}(t)=1/V\mathcal{F}(it)=\mathcal{F}(-t)=\mathcal{F}(t)$.
\begin{proposition}
\label{prop:invol}
 The mapping $V$ defined by \eqref{eq:defD} is an involution on $\D$.
\end{proposition}
In other words,   the set $\Dt$ is closed under commutation, that is, if $[\mathcal{F},V \mathcal{F}] \in \Dt$
then $[V \mathcal{F},\mathcal{F}] \in \Dt$.

Next, from the Definition \ref{def:DL} of $\De$, a first natural question is to  understand
whether the  Lukacs class $\De$ contains  all possible  entire characteristic  functions with
only real zeros, that is, whether there exists such an entire function of order $\rho>2$. Here is the definitive answer to this issue,  which is, in fact, a direct consequence of a very nice  result due to Gol'dberg and Ostrovs'ki  \cite{Goldo}, see also
{\cite[Theorem 4.4.1]{Lukacs_book}},
%\tkside{Supply reference}
regarding entire characteristic functions having only real zeros. Note that, in these references,  the statement is proven for a more general class of entire functions, namely the ones  possessing the so-called ridge property, that is,
$|\varphi(z)| \le |\varphi(\Im(z))|$.

\begin{proposition}
\label{DLD}
Every entire characteristic function with only real zeros belongs to the class
$\mathds D_L$.
\end{proposition}
%\begin{theorem}
 %$\mathcal{F}$ is an entire function with only real zeros in $\Pp_+ $ if and only if $\mathcal{F} \in \De=\LPnew \cap \Pp_+$.
%\end{theorem}
%\begin{proof}
%  Let $\mathcal{F}$ be an entire function with only real zeros in $\Pp_+ $. First, according to   the order $\rho$ of the entire function $\mathcal{F}$ is no greater than $2$. This yields that $\mathcal{F} \in \mathds{LP} \cap \Pp_+$. According to Schoenberg \cite[Theorem 1]{Schoenberg}, $\frac{1}{\mathcal{F}(it)} \in \Pp_+ $
%\end{proof}
We proceed by deriving some closure properties of the sets $\D$ and $\De$.
First, since  multiplication of characteristic functions remain characteristic function as they
correspond to addition of  independent random variables, we get that $\D$ is stable by multiplication and having the constant function $1$ as identity element, we obtain that it is a monoid.
Similarly, since any reciprocal of a function in $\De$ being  the moment generating function, at least on an imaginary strip, of a P\'olya frequency function,  it remains to identify transforms that preserve the positive definiteness property of a Laguerre-P\'olya function. It is the program that we develop in the remaining part of this section.  In this spirit, there is this first closure property which follows readily since the product of two characteristic functions in $\mathds{LP}$, that is of the form \eqref {eq:defPL}, remain an even entire  characteristic function in $\mathds{LP}$, see the item \eqref{it:mono} above for the same property for the set $\D$.
%\tkside{Alternatively, if $F(z)$, $G(z)$ are\mbreak even characteristic functions
%\mbreak in $\mathds{LP}$ then $F(z)G(z)$ is also \mbreak
%an even ch.f.\ and in $\mathds{LP}$.}
\begin{proposition}\label{it:mono}
The sets $\D$ and $\De$ constitute a monoid under the operation of function multiplication.
\end{proposition}
Note that the previous claim could also be interpreted as  the set $\De$ being invariant by the convolution of probability distributions or equivalently by taking the sum of independent random variables. The following mappings and results  were proposed by Lukacs
\begin{equation}\label{eq:lukacs}
L^{(p)} f(t)= t^{p-2}\frac{f^{(p)}(t)}{f^{(2)}(0)}, \quad t\in \R,\quad
p=1,2,
\end{equation}
where, for a
 sufficiently smooth
function $f$, we write $f^{(p)}(t)=\frac{d^p}{dt^p}f(t)$.
\begin{theorem}\cite[Theorem 4]{Lukacs}\label{it:dif}
We have, for $p=1,2$,   $L^{(p)}( \mathds{LP})\subset \mathds{LP}$
and $L^{(p)}( \De)\subset \De$.
%   where $L^{(p)} f(t)= t^{p-2}\frac{f^{(p)}(t)}{f^{(2)}(0)}, t\in \R$,  where, for a smooth function $f$, we write $f^{(p)}(t)=\frac{d^p}{dt^p}f(t)$.
\end{theorem}
The proof of this theorem relies on Laguerre's theory, see
e.g.\ Borel \cite[Ch.\ 2]{Borel}, which provides a  closure property of the class
$\mathds{LP}$ by differentiation.
Then Lukacs showed, by analytical means,  that his mappings leave the set
$\Pp_+$  invariant.  In Section \ref{sec:proof} we shall provide
an alternative proof of this last fact based on probabilistic arguments.

\begin{remark}
\begin{enumerate}[(i)]
  \item Since $L^{(p)}$, $p=1,2$, are differential operators, we easily get that the unique invariant of $L^{(1)}$ (resp.~$L^{(2)}$) in $\Pp_+$, i.e.~$L^{(1)} \mathcal F = \mathcal F$, is
$\mathcal F(t) = e^{-\sigma^2 t^2/2}$ (resp.~$\mathcal F(t) = \cos(\sigma t)$, $\sigma\in \R$).
\item For any finite sequence $(p_1, \ldots, p_n)$ whose elements take values in
$\{1,2\}$ we have that
 $L^{(p_n)} \circ \cdots \circ L^{(p_1)} (\D_L)\subseteq \D_L$. For example, one observes that
\[
L^{(1)}\cos(t) =\frac{\sin t }{t} \textrm{ and }
L^{(1)} \circ L^{(1)}\cos(t) = 3 \frac{\sin t - t \cos t}{{t}^3}.
\]
More generally, define $\mathcal F_{n+1}(t) = L^{(1)} \mathcal F_n(t), n\geq0$, with $\mathcal F_0(t)=\cos t$.
Thus, $\mathcal F_n$ is obtained by the $n$-fold application of the operator $L^{(1)}$ to
the cosine function.
We obtain the expression
\[
\mathcal F_{n+1}(t) = \frac{-(2n+1)!!}{t^{2n+1}}\,  \Re (\overline{P}_n(-it) i e^{-it})
\]
where $(2n+1)!! = (2n+1) (2n-1) (2n-3) \cdots 1$
and $\overline{P}_n$ is the so-called reverse Bessel polynomial,
expressed in terms of the degree-$n$ Bessel polynomial $P_n$ via
$\overline{P}_n(x) = x^n P_n(1/x)$, see \cite{KF48}.
Note that $\mathcal F_n \to \mathcal F_0$ pointwise as $n \to \infty$ and this is natural
from the probabilistic interpretation of the operator $L^{(1)}$, see Lemma \ref{lem:Lm}.

\end{enumerate}
\end{remark}

We proceed with the following two results that give a characterization or a partial characterization of the set $\De$.
First, we have the following  L\'evy-Khintchine type representation which follows readily from Theorem \ref{thm:sch} combined with the characterization  result due to Kwa\'snicki  of the Laplace transform of P\'olya frequency densities.
\begin{theorem} \cite[Proposition 5.3]{Mateu}
If $\mathcal{F} \in \De$ then, for all $t\geq0$,
\label{eq:polya:exp}
\begin{equation}
 \mathcal{F}(t)  = e^{\Psi_L(t)}
\end{equation}
where  $\Psi_L(t) = -c t^2 - \int_{-\infty}^\infty \left( \frac{1}{t + r} - \frac{1}{r} + \frac{t}{r^2} \right) \rho(r) dr, c \ge 0$ and $\rho : \R \to \Z$ is an even non-decreasing integer-valued function such that $
 \int_{-\infty}^\infty \frac{|\rho(r)|}{|r|^3} \, dr < \infty$.
\end{theorem}

The next theorem, proved by Newman \cite{New1},  who was, as discussed in Section
\ref{sec:LY}, motivated by problems arising in statistical physics and
Euclidean field theory, characterizes a subclass of $\mathds D_L$.
For a random variable $\mathrm X$ and a real number $\lambda$ such that
$Z_\lambda := \int_{-\infty}^\infty e^{\lambda x^2} F_{\mathrm X}(dx) < \infty$,
let $\mathrm X_\lambda$ denote a random variable whose distribution
$F_{\mathrm X_\lambda}$
is given by
\begin{equation}\label{eq:defF}
F_{{\rm X}_{\lambda}}(dx)
= \frac{1}{Z_\lambda} e^{\lambda x^2} F_{\rm X}(dx), \: x\in \R.
\end{equation}

\begin{theorem}\cite[Theorem 1]{New1}
\label{thm:New}
Let $X$ be a symmetric random variable. Then,  $\mathcal F_{{\rm X}_{\lambda}} \in \De$ for all $\lambda \in D_{\rm X}=\{\lambda \in \R; \: \mathcal F_{\rm X^2}(-i\lambda)<\infty  \}\supseteq (-\infty,0]$ if and only if either, for some $x_0$,
\[
F_{\rm X} =\frac{1}{2}\left( \delta_{x_0} +\delta_{-x_0} \right) ,
\]
or  $F_{\rm X}$  is absolutely continuous with respect to the Lebesgue measure with a density $f_{\rm X}$ which takes the form
\begin{equation}
f_{\rm X}(x) = Kx^{2m}e^{-\alpha x^{4}-\beta x^{2}} \prod_{k=1}^{N}\left( \left( 1+%
\frac{x^{2}}{a_{k}^{2}}\right) e^{-\frac{x^{2}}{a_{k}^{2}}}\right)
\label{DBNdensity}
\end{equation}%
where $K>0$ is a normalizing constant, $m$ a nonnegative integer,
$\alpha, \beta$ real numbers, $N$ a nonnegative integer or $\infty$,
with the $a_k$ positive, and
either $\alpha=0$, $\sum a_k^{-4}<\infty$,
or $\alpha >0$, $\beta + \sum a_k^{-2}>0$ (the case $\sum a_k^{-2}=\infty$
is allowed).
\end{theorem}

\begin{remark}
This theorem does not fully characterize the set $\De$ as, for instance,
$\mathcal{J}_{0}$, the Bessel function of order $0$, see \eqref{eq:ex2}, is in  $\De$,
but it is the characteristic function of the arc-sine  law whose density, see \eqref{eq:defJ}, does not have the form \eqref{DBNdensity}. Another
 well-known function that does not belong to Newman characterization  is the function
\[
t\mapsto\mathcal{F}_{\lambda,\Phi}(t)= \frac{1}{\int_{\R}e^{\lambda x^2} \Phi(x)dx}\int_{\R}  e^{itx}e^{\lambda x^2} \Phi(x)dx,
\]
introduced by P\'olya,  where $\Phi$ is the inverse Fourier transform
of the Landau function $\xi$; see \eqref{xitrans}.
It is well-known that \[\mathcal{F}_{\lambda,\Phi} \in \De \textrm{ if and only if } \lambda\geq\Lambda_{DN}\] where $\Lambda_{DN}$ is the celebrated de Bruijn-Newman constant. Since,
by Theorem \ref{thm:RH},
the Riemann hypothesis is equivalent to $\Lambda_{DN}\leq 0$, this observation has motivated an intensive research activity on the computation of $\Lambda_{DN}$. The current  state of art is $ 0 \leq \Lambda_{DN}\leq  0.22 $ and was obtained by Tao and collaborators \cite{RT,Tao}.
\end{remark}

The following three claims are adaptation to our setting of some deep results  due to P\'olya \cite{Polya27} (the first two) and to de Bruijn \cite{DB}, see also \cite[sections 2.2 and 2.3]{NewRH}. We omit their proofs as they follow readily from the aforementioned results combined with  Theorem \ref{thm:sch}.

\begin{theorem}\cite{Polya27}\label{thm:p1}
Let $f: \R \to \R_+$ be the density of a symmetric random variable such that,
for some $A, \alpha >0$,
\begin{equation}\label{eq:buf}
f(x) \leq A e^{-x^{2+\alpha}}, \: x \ge 0.
\end{equation}
Let $\varphi : \R \to \R_+$ be even, real analytic, and such that
$\int_{\R}  f(x)\varphi(x)dx=1$.
Assume further that $t\mapsto\mathcal{F}_f(t) = \int_{\R}  e^{itx} f(x)dx \in \De$.
Then,
\[
t \mapsto \mathcal{F}_{f\cdot \varphi}(t)=\int_{\R}  e^{itx} f(x)\varphi(x)dx\in \De
\]
if and only if the analytic extension of  $\varphi$ is such that  $t\mapsto \varphi(it) \in \mathds{LP}$.
\end{theorem}

Note that the function $\varphi$ in the above theorem is called a universal factor by
P\'olya. An interesting and simple application of this theorem is for the entire function
$\varphi(x)=e^{\lambda x^2}$ where $\lambda >0$.
Indeed, using the notation of Theorem \ref{thm:New},
if ${\rm X}$ is a symmetric random variable with an absolutely continuous distribution
whose density satisfies the bound \eqref{eq:buf} and  for some
$\underline{\lambda} \in \R$, $\mathcal F_{{\rm X}_{\underline{\lambda}}} \in \De$,
then $\mathcal F_{{\rm X}_{\lambda}} \in \De$ for all
$\lambda\geq \underline{\lambda}$ as plainly here $D_{\rm X}=\R$.
For instance, from the expression \eqref{eq:defJ},
we deduce that, for any $\nu<\frac12$ and $\lambda>0$, the mapping
\begin{equation}
t \mapsto K_{\lambda} \int_{|x|<2}
e^{itx}e^{\lambda x^2}(4-x^2)^{{-\nu}-\frac{1}{2}}dx \in \De
\end{equation}
where $K_{\lambda}>0$ is a normalizing constant.
However, note that since the density above when $\lambda=0$ is not of the form
\eqref{DBNdensity}, by Theorem \ref{thm:New},
there exists $\underline{\lambda} <0$ such that the entire function
\begin{equation}
t \mapsto K_{\underline{\lambda}} \int_{|x|<2} e^{itx}
e^{\underline{\lambda} x^2}(4-x^2)^{{-\nu}-\frac{1}{2}}dx
\end{equation}
has non-real zeros.

\begin{theorem}\cite{DB}      \label{thm:p2}		
Let $f: [0,\infty) \to \R_+$ be the density of a positive random variable,  %hermitian, $f(-t)=\overline f(t)$,
and, for all $y\ge 0$,
\[
f(y) \leq B e^{-y^{\frac12+\beta}}
\]
for some $B, \beta >0$.
Suppose that $f$ has an analytic extension in a neighborhood of the origin.
Then  the function
\[
\mathcal{M}_f(t) = \int_0^\infty y^{t-1} f(y) dy
\]
has an extension on $\C$ as a meromorphic function.
If the function $\mathcal{M}_f$ has only negative zeros and $n$ is a positive integer,
then, writing $\frac{1}{C_n}=\int_{\R}   f(x^{2n})dx$,
\[
t \mapsto C_n\int_{\R}  e^{itx} f(x^{2n})dx \in \De.
\]
%is an entire function with only real zeros.
\end{theorem}						

\begin{theorem}\cite{NewRH}			\label{thm:p3}
Let $f$  be an entire function such that its derivative $f^{(1)}$ is the limit
(uniform in any bounded domain) of a sequence of polynomials, all of whose roots
lie on the imaginary axis. Suppose further that $f$ is not a constant, $f(x) = f(-x)$,
and $f(x) \geq0$  for $x\in \R$ and $\int_{\R}
e^{-f(x)}  dx=1$. Then
\[ t\mapsto \int_{\R}e^{itx}
e^{-f(x)}  dx \in \De.\]
\end{theorem}

As instances illustrating these results, there are the following entire functions that were derived by P\'olya
\begin{eqnarray}
 t\mapsto  &K\int_{\R}e^{itx}\cosh (ax) e^{ -a\cosh x} \,dx, &
 t\mapsto K\int_{\R}e^{itx}e^{ -x^{2n}} \,dx \textrm{ and } \nonumber \\
 t\mapsto &K\int_{\R } e^{itx}e^{-ax^{4n}+bx^{2n}+cx^{2}}
\,dx \in \De& \label{Polyaexample}
\end{eqnarray}%
where $K$ is, in each expression, a normalizing constant, and $n\in \mathbb{N}, a>0$.% and so does the entire function
%\[
%\begin{equation}
%\int_{-\infty }^{\infty } e^{itx}\exp\left( -ax^{4q}+bx^{2q}+cx^{2}\right)
%\,dx,\text{ }q\in \mathbb{N}\text{.}
%\end{equation}

We continue P\'olya's and de Bruijn's line of research
by presenting original additional closure properties of the set $\De$. More specifically, we
%introduce the following notion of random mixing of elements of $\De$, which can be seen
%as
investigate its stability  under product of independent variables.
We  have
already mentioned that this property is intimately connected to the concept of intertwining
relationship between Markov semigroups, as we will discuss later in the paper.
To state it, we say  that  a positive linear operator $\Lambda$
 on the space of bounded borelian functions
%\tkside{In other words, if we let\mbreak ${\rm Y}_n = {\rm I}_0 {\rm I}_1\cdots {\rm I}_n$ then $\Lambda_{\rm I}$ is the \mbreak generator of the semigroup \mbreak of the Markov chain $({\rm Y}_n)$.}
is Markov multiplicative if there exists a random variable ${\rm{I}}$ with distribution
function $F_{{\rm{I}}}$, such that, for any bounded borelian function $f$, writing
$\Lambda=\Lambda_{{\rm{I}}}$,
\begin{equation}\label{eq:defLI}
\Lambda_{{\rm{I}}} f(t)=\int_{\R}f(xt)F_{{\rm{I}}}(dx).
\end{equation}

In the following, we  identify a mapping from the set of even entire functions in $\Pp_+$ into $\De$. In other words,  we provide a way of creating entire characteristic functions with only real zeros from any even characteristic functions. We also find necessary conditions on a Markov multiplicative operator to leave invariant the set $\D_L$.
Let us now denote by $\mathds{LP}_+$ the class of functions in  $\mathds{LP}$
with only strictly negative zeros and of the form
\[
\varphi(z) = K e^{az} \prod_{k=1}^\infty
\left(1+\frac{z}{z_k}\right) e^{-z/z_k}
\]
with $K,  a \in \R$ and $z_k > 0$ for all $k$ with
$\sum _{k \ge 1} z_k^{-2} < \infty$.

\begin{theorem} \label{thm:DM}
\begin{enumerate}
 \item \label{it:Lag1}Let  $\mathcal{F}_{\rm{D}}  \in \Pp_+$ be entire and even, and,   assume that there exist
 a random variable $\rm{I}$ such that for any non-negative integer $n$
\[
\mathcal{M}_{\rm{I}}(2n)\mathcal{M}_{\rm{D}}(2n) =a_{\varphi}(n)G(n)\Gamma(2n+1)
\]
where {$\varphi(z)=\sum_{n=0}^{\infty}a_{\varphi}(n)z^n \in \mathds{LP}$ and
$G \in \mathds{LP}_+$}.
{ Then,} $\Lambda_{{\rm{I}}}\mathcal{F}_{\rm{D}} \in \De$.
\item \label{it:Lag2} Let us now assume that there exist      $\varphi \in \mathds{LP}_+$ and  a random variable $\rm{I}_L$ such that for any non-negative integer $n$,
\[
\varphi(n)=\mathcal{M}_{\rm{I}_L}(2n).
\]
{ Then,}  $\Lambda_{{\rm{I}}_L}(\De) \subset \De$. %, that is, $\Lambda_{{\rm{I}}}$ is a van Dantzig mixing operator that preserves the class $\De$.
    \end{enumerate}
\end{theorem}
Before continuing further, let us first illustrate
 part (1) of this theorem
by  an example. Let us  take
$\varphi(z)=e^{z}\in \mathds{LP}$
%\tkside{Should this be $\mathds{LP}_1$?}
and then $a_{\varphi}(n)=\frac{1}{n!}, n\geq0$.  For instance,  let
${\rm{D}}={\rm{J}}_0$ where ${\rm{J}}_0$ is the symmetric random variable whose distribution is the arc-sine law and  is recalled in \eqref{eq:defJ} below. Then, from \eqref{eq:MD}, $
\mathcal{M}_{{\rm{J}}_0}(2n)=\frac{(2n)!}{n!n!}$ and $\mathcal{F}_{{\rm{J}}
_0}= \mathcal{J}_{0} \in \De$, the Bessel function of order $0$, see \eqref{eq:ex2} above.

 Next fix $b>0$ and define ${\rm{I}}_b$ as the positive
random variable
 with distribution
$F_{\rm{I}_b}(dx)=\frac{2b}{b+1} e^{-
x^2}x\left(\frac{x^2}{b}+1\right)dx, x>0$. Simple
algebra yields that $\mathcal{M}_{\rm{I}_b}
(2n)=\frac{(n+b)}{b}n!$, and, with the previous choice of $\varphi$,
the equation
$\mathcal{M}_{\rm{I_b}}(2n)\mathcal{M}_{\rm{J_0}}(2n) =a_{\varphi}(n)G(n)\Gamma(2n+1)$
gives
\[
G(n)=\frac{(2n)!n!(n+b)n!}{bn!n!(2n)!}= \frac{(n+b)}{b}.
\] Since $G(z)=\frac{z+b}{b}$ is in  $\mathds{LP}_+$,
the theorem above shows that

\[
\Lambda_{{\rm{I}}_b}\mathcal{F}_{{\rm{J}}_0}(t)=\sum_{n=0}^{\infty}(-1)^n\frac{n+b}{n!b} t^{2n}=\frac{b-t^2}{b}e^{-t^2}\in \De,
\]
for all $b>0$. Note that this characteristic function already appeared in  \cite[Example III.8]{Linnik77} where the authors showed that it is not decomposable with respect to the additive convolution of probability measures. However, our approach reveals that it is decomposable with respect to the multiplicative one.

We point out that Theorem \ref{thm:DM} enables one to generate many new examples of elements in $\De$, and, we refer to Theorem \ref{thm:patie} where this idea is exploited and to Section \ref{sec:exa} below, where some additional examples are provided.

Another interesting aspect of the set $\De$ is the following additional closure property.
\begin{theorem} \label{thm:conv}
The sets $\D$ and $\De$  %=\{\mathcal{F}_X \in \De; \mathcal{F}_{X^2}(-iu)<\infty \textrm{ for some } u>0 \}$
are closed under locally uniform  convergence.
\end{theorem}
Note that  the locally uniform  convergence  of characteristic functions is equivalent to   the pointwise convergence to a continuous function  at $0$, see \cite[Theorem 3.2.1]{Boch}, and, by the L\'evy continuity theorem, see \cite[Theorem 8.28]{Brei}, this implies  the weak convergence of the corresponding sequence of random variables to a random variable. In other words, Theorem \ref{thm:LY} entails that Theorem \ref{thm:conv} is a generalization  to \cite[Theorem 7]{New} regarding the closure under weak convergence of the set $\mathcal{X}$, which is the set of symmetric probability measures $F$ whose characteristic function is in $\LPnew$ and such that $\int_\R e^{bx^2}F(dx)<\infty$ for some $b>0$. We do not need this last condition for our closure result. However, the result from Newman and Wu \cite[Theorem 7]{New} ensures that the gaussian tail property  for elements in $\De$ is preserved under weak convergence.

\section{The new class $\D_P$}
\label{intro part 2}
To present the main results of this Section, we start by introducing  some objects and notation.
First, with  $\mathbb{M}_+(\R^+)$ denoting the set of non-negative Radon measures
on $(0,\infty)$, we define the mapping $\Psi :\R^+\to \R$ by
\begin{equation}
\Psi(u) =-\kappa+au +\frac{1}{2}\sigma^{2}u^{2}
-\int_0^{\infty}(1-\mathrm{e}^{-u r}-u r\mathbb{I}_{\{r<1\}})\mu ({\rm d}r)
\label{LLK}
\end{equation}
where $\kappa\geq 0$, $a\in\mathbb{R}$, $\sigma^2\geq 0$ and
$\mu\in \mathbb{M}_+(\R^+)$ is such that
$\int_{0}^{\infty} \inf(1, r^2) \mu({\rm d}r)<\infty$.
We exclude the case when $\sigma=0$,
$a  + \int_0^{1} r\mu ({\rm d}r)\leq0$ and
$\int_{0}^{\infty} \inf(1, r) \mu({\rm d}r)<\infty$,
 which is seen as degenerate in our context.
This function has a nice probabilistic interpretation.
Indeed, it is the so-called Laplace exponent of a possibly killed real-valued
spectrally negative L\'evy process $Y=(Y_t)_{t\geq 0}$,
i.e.~a stochastic process without positive jumps and with stationary and independent
increments starting from $0$ and, when $\kappa>0$,
it is killed at an independent exponential time of parameter $\kappa$.
Moreover, we have, for any $u,t\geq0$,

   \begin{equation}\label{eq:LeY}
     \mathcal{F}_{Y_t}(-iu)=e^{\Psi(u)t}.
   \end{equation}

  Note that under the three conditions we excluded above, $Y$ is negative-valued and has non-increasing sample paths. We refer to the monograph \cite{Kyprianou06} for a thorough study of these processes.
It is a well established (and an easy to check) fact that the Laplace exponent $\Psi$ is  strictly convex on $[0,\infty)$ with
\begin{equation}\label{eq:propPsi}
  \lim_{u\to\infty}\Psi(u)=+\infty \textrm{ and } \Psi \textrm{ is increasing on } [\theta,\infty)\textrm{ where }
 \theta=\sup\{u\geq0;\, \Psi(u)=0\}.
\end{equation}
In fact, by convexity, we have
\begin{equation}\label{eq:zeta}
\textrm{$\theta>0$  if and only if
(i) $\kappa>0$ or
(ii) $\kappa=0$ and $\Psi^{(1)}(0^+)=\lim_{u \downarrow 0 }\Psi^{(1)}(u) <0$}
\end{equation}
and a monotone convergence argument yields that
$\Psi^{(1)}(0^+)={a {-} \int_1^{\infty}r\mu ({\rm d}r)}
\in [-\infty,+\infty)$.
We are now ready to define the sets
\begin{equation}\label{eq:defNG}
\N=\{ \Psi :\R^+\to \R \textrm{ of the form } \eqref{LLK} \}
\end{equation}
and
\begin{equation}\label{eq:defN}
\N_{\D}=\left\{ \Psi \in \N;\, 0\leq \theta {{ \leq}} \tfrac12 \right\}.
\end{equation}
Note that  \eqref{eq:zeta} entails that
\begin{equation}\label{eq:propNd}
\Psi \in \N_{\D}  \textrm{ if } \kappa=0  \textrm{ and } \Psi^{(1)}(0^+)\geq0 \end{equation}
as, in this case, $\theta=0$. Moreover, from \eqref{eq:propPsi}, one easily gets that a necessary and sufficient condition for a function $\Psi$ to be in   $\N_{\D}$ is that $\Psi \in \N$ with $\Psi(\frac12){{ \geq }} 0$.
The notation of the set $\N_{\D}$ is motivated by the following facts.
On the one hand, due to the infinite divisibility of $Y_1$,  we have, for any $\Psi \in \N$, that
$
z\mapsto -\Psi(iz) \in \N(\mathcal{\R}),
$
the set of continuous and negative-definite functions on $\R$.
On the other hand, we now define a class of entire functions that are generated by the set
$\N_{\D}$, that will be shown to  belong to $\D$.

For any $\Psi \in \N$, we introduce the power series
%\footnote{ E.g., with $\Psi(u)=u^2/2$, corresponding to standard BM,
%we get $\mathcal J_\Psi(t) = J_0(t)$, the first-kind Bessel function of order $0$.}
\begin{equation}\label{eq:defIp}
 \mathcal{J}_{\Psi}(t)=\sum^{\infty}_{n=0} \frac{(-1)^{n}  }{W_{\Psi}(n+1)}t^{2n}
\end{equation}
where, here and below, for a  function $\varphi$ defined on $\R^+$, we  set
\[
W_{\varphi}(1)=1, \quad W_{\varphi}(n+1)=\prod^{n}_{k=1} \varphi(k),
\quad n \ge 1.
\]
Since, from \eqref{eq:propPsi}, $\lim\limits_{n\to\infty}\frac{W_{\Psi}(n+2)}{W_{\Psi}(n+1)}=\lim\limits_{n\to\infty}\Psi(n+1)=\infty$, $\mathcal{J}_{\Psi}$ defines an entire function.

We also write
\begin{equation}
\mathcal{I}_{\Psi}(t)=\mathcal{J}_{\Psi}({ it}
%e^{i\frac{\pi}{2}}t
)=\sum^{\infty}_{n=0} \frac{1  }{W_{\Psi}(n+1)}t^{2n}
\end{equation}
 and point out that the entire function $\mathcal{I}_{\Psi}(\sqrt{t})$ was introduced by the second author in \cite{Patie} where it was shown that it is an  invariant function of some self-similar integro-differential operator.
Therein,  the complete monotonicity property  was identified for several of its transformations,  and, we also refer to \cite{Bart-Patie} for a more recent and refined studies of this class of functions. We are now ready to introduce the following class of entire functions
\begin{equation}\label{eq:defDp}
  \Dp=\{\mathcal{J}_{\Psi} \textrm{ of the form } \eqref{eq:defIp} \textrm{ with } \Psi \in \Nd\}.
\end{equation}

\begin{theorem} \label{thm:main}
We have $\Dp \subseteq \D$, and, for all $\Psi \in \N_{\D}$,
\begin{equation}
\left[\mathcal{J}_{\Psi},\frac{1}{\mathcal{I}_{\Psi}}\right]  \in \D_2
\end{equation}
with
\begin{equation}
t\mapsto - \log \mathcal{J}_{\Psi}(t) \notin \N(\R)
\textrm{ but } t\mapsto   \log \mathcal{I}_{\Psi}(t) =\phi_{\Psi}(t^2) \in \N(\R)
\end{equation}
where the function $\phi_{\Psi}$ is a Bernstein function that belongs to the class $\Be_J$,
that is, it is of the form  \eqref{eq:bernJ} below with the additional property
that $r\mapsto r \overline{\mu}(r)$ is non-increasing on $\R^+$. %,   such that, for all $t\in\R$,
Finally, writing $ \mathcal{F}_{\Psi}(t)=\frac{\mathcal{J}_{\Psi}(t)}{\mathcal{I}_{\Psi}(t)}, t\in \R$, we have  $\mathcal{F}_{\Psi} \in  \D_S$, that is $\mathcal{F}_{\Psi}(t)\mathcal{F}_{\Psi}(it)=1$ for all $t\in\R$.
\end{theorem}

  \begin{remark}
%We shall describe in the proofs below the random variables whose characteristic function appeared in the previous result.
 The random variables whose characteristic functions appear above shall be
explicitly described in Section \ref{charprop}.
  \end{remark}

  \begin{remark}
Note that in \cite{PKS}, it is proved that the entire function
$\mathcal{J}_{\Psi}(\sqrt{t})$ has its smallest (in modulus) zero, say $z_1$,
which is simple and located on the positive real line.
On the one hand, this shows that $\mathcal{J}_{\Psi}$ is not the
characteristic function of an infinitely divisible variable
as their characteristic functions are zero-free,
see \cite{Sato}. On the other hand, $z_1<0$   corresponds to
a singularity of the Bernstein function $\phi_{\Psi}$.
However, from Theorem \ref{thm:main}, we get the identity, for all $t\in \R$,
\begin{equation}\label{eq:prodvd}
   \mathcal{I}_{\Psi}(t) \mathcal{J}_{\Psi}(t) =e^{-\phi_{\Psi}(t^2)}e^{\phi_{\Psi}(-t^2)}=1,
\end{equation}
which  entails that the Bernstein function
$\phi_{\Psi}$ admits a meromorphic extension on $\C$, but it is not necessarily a Pick function (see the definition below). This reveals that such an identity cannot be possible for  Bernstein functions, especially for those having an essential singularity, e.g.~$\phi(u)=u^{a}, 0<a<1$.
  \end{remark}

The purpose of the next result is to explain how the classes  $\Dp$ and $\De$ are related, which consists on investigating  the difficult issue of locating  the zeros of the entire function $\mathcal{J}_{\Psi} \in \Dp$.
 To state it, we recall that a Bernstein Pick function
\cite[p.\ 56]{SchillingSongVondracek10}
is a Bernstein function which admits an holomorphic extension which maps the upper half-plane into its closure. We say that an entire function  (resp.~Pick meromorphic function)  has the $1$-separation property if its sequence of zeros $(z_k)_{k\geq 1}$ (resp.~and poles $(\rho_k)_{k\geq 1}$) satisfies, for all $k$, $z_{k+1}< z_k-1$ (resp.~$\rho_{k}=z_{k}-1>z_{k+1}$). We now  introduce the set
\begin{equation}
\Be_{P_1}=\{\phi \in \Be; \phi \mbox{ is a  Pick function having the $1$ separation property}\}.
\end{equation}
\begin{theorem} \label{thm:pick}
Let $\phi\in\Be_{P_1}$. Then, $u\mapsto\Psi(u)=u\phi(u) \in \Nd$ and $\mathcal{J}_{\Psi} \in \De \cap \Dp$. However, $\Dp \nsubseteq \De$ as there are $\Psi's \in \Nd$ such that $\mathcal{J}_{\Psi}$ has at least a non-real zero.  %Moreover,  the Riemann $\xi$ function does not belong to $\Dp$.
\end{theorem}
\begin{remark}
  In Section \ref{sec:exa}, we provide instances of the two situations presented in this theorem, see e.g.~ the Bernstein functions that define the Bessel functions and the Fox-Wright functions, see \ref{sec:fw}, which both belong to  $\De$. On the other hand there are the examples involving hypergeometric functions and the Mittag-Leffler functions which have
 non-real zeros.
\end{remark}

 We proceed with the following
 result
that shows that the Lukacs mappings, introduced in \eqref{eq:lukacs},
also leave our class $\Dp$ invariant.
  \begin{proposition}\label{prop:LukMap}
 Let  $L^{(p)}, p=1,2,$ be the operators that were defined in \eqref{eq:lukacs} above. Then, we have  $L^{(1)}(\Dp)\subset \Dp$. Moreover, $L^{(2)}(\Dp)\subset \D$ and the same  remains true for their iterates.
%More specifically, we have $T^{p}\mathcal{J}_{\Psi}(t)=\frac{1}{W_{\Psi}(p+1)}((-2t)^{-p}\mathcal{J}_{\Psi}(t))^{(p)}$.
  \end{proposition}

\subsection{The Riemann $\xi$ function and the class $\Dp$}\label{sec:vDR}
A natural and important question that arises at this stage is to understand whether the Riemann $\xi$ function defined in \eqref{xitrans} belongs to the class $\Dp$. Indeed, this would yield  a power series representation of this function whose coefficients would be expressed in terms of negative definite functions offering new tools to study the location of its zeros, using for instance Theorem \ref{thm:pick}. To this end, let us recall that the Riemann $\xi$ function, defined in \eqref{xitrans}, can be expressed in terms of the following power series
\begin{equation} \label{eq:xi}
	\Theta(z)=\xi(\sqrt{z})=\sum^{\infty}_{n=0} \frac{\gamma(n)}{n!}z^{n}, \ z\in \mathbb{C},
	\end{equation}
where, with  $ F(n) =  \int_1^\infty (\log x)^n\,x^{-3/4}\,\theta_0(x)\,dx$ and $\theta_0(x)=\sum_{n=1}^{\infty} e^{-\pi n^2x}$  the theta series, we have set
\begin{equation}\label{gamma_formula}
\gamma(n)=\frac{n!}{(2n)!} \frac{32 \binom {2n}2F(2n-2) - F(2n)}{2^{2n-1}},
\end{equation}
see e.g.~\cite{Ono}.
The question whether there exists $\Psi \in \Nd$ such that $\mathcal{J}_{\Psi}=\xi$ boils down to the existence of $\phi \in \Bd$ such that, for all $n\in \mathbb{N}$,
\begin{equation}\label{gamma_formula}
\frac{1}{W_\phi(n+1)}=\frac{n!}{(2n)!} \cdot\,\frac{32 \binom {2n}2F(2n-2) - F(2n)}{2^{2n-1}},
\end{equation}
which, after some easy algebra, is equivalent to show that
\begin{equation}\label{gamma_formula}
\phi(n+1)=-\frac{G(2n)}{8(n+1)G(2n+2)},
\end{equation}
where we have set $ G(2n)= 64n(2n-1)\frac{F(2n-2)}{F(2n)}-1$.
 Since this question does not seem straightforward,  we investigate, instead here, whether this possibility could be excluded from the properties that we know about  elements of the class $\Dp$ and the Riemann $\xi$ function. Recalling that the latter is an entire function of order $1$ with infinite type and it is the characteristic function of the density of a  probability measure whose support is  $\R$, we have the following.
\begin{proposition}
Let $\Psi \in \Nd$ such that $\Psi(u)=\frac{u^2}{\ell(u)}, u\in \R,$  with $\lim_{u\to\infty}\ell(u)=\infty$, $\ell$ being a slowly varying function at infinity, i.e.~for every $u>0$, $\lim_{t\to\infty}\frac{\ell(ut)}{\ell(t)}=1$.  Then $\mathcal{J}_{\Psi}$ is an entire function of order $1$ and infinite type. This condition holds when $\Psi(u)=(u-\theta)\phi(u), \theta\geq0,$ with $\phi$ a special Bernstein function, i.e.~$\phi(u)\widehat{\phi}(u)=u$, such that its conjugate Bernstein function $\widehat{\phi}(u)=\ell(u)$.   Moreover, under this condition,  $\mathcal{J}_{\Psi}$ is the characteristic function of a probability density function whose support is  $\R$. %, is its so-called Blumenthal-Getoor lower index.
\end{proposition}
The first part of the Proposition follows readily from \cite[Proposition 2.1]{Bart-Patie}, see also Proposition \ref{prop:order} below, where we notice that when the order is $1$, with the notation of  Proposition \ref{prop:order}, $\underline{\Psi}=2$ and thus the type $\tau_{\Psi}\geq \left(\limsup_{n \to \infty}\frac{n}{\phi(n)}\right)^{\frac{1}{2} } =\left(\limsup_{n \to \infty}\ell(n)\right)^{\frac{1}{2}}=\infty$ under the condition of the Proposition. The last claim is a specific instance of Lemma \ref{lem:IJ}.
\begin{remark}
One instance when $\Psi(u)=\frac{u^2}{\ell(u)}, u\in \R,$ with $\ell$ as in the Proposition is when $\Psi(u)=ue^{{\rm{W}}(u)}$, where ${\rm{W}}$ is the Lambert function. Indeed, it is well known that for all $u \geq0$, ${\rm{W}}(u)e^{{\rm{W}}(u)}=u$ and ${\rm{W}}$ is a complete Bernstein function that is a Bernstein function whose L\'evy measure is absolutely continuous  with a completely monotone  density, and  $\lim_{u\to \infty}\frac{{\rm{W}}(u)}{\ln u}=1$, and hence it is a special Bernstein function with  $e^{W(u)}$ as conjugate, see \cite{Pakes}. Therefore  $\Psi(u)=ue^{{\rm{W}}(u)}\sim\frac{u^2}{\ln u}$.
\end{remark}

\subsection{Some  examples in the class $\D_P$}\label{sec:exa}
In this section, we give several specific examples of the function $\mathcal{J}_{\Psi}$ including the modified Bessel functions, the Mittag-Leffler functions and several type of hypergeometric functions, and refer to \cite{Gasper, Kilbas06} as classical references on these functions.  The interested reader can also consult  the monograph \cite{SchillingSongVondracek10} for several examples of Bernstein functions from which one can provide additional interesting instances of $ \mathcal{J}_{\Psi}$.

\subsubsection{Bessel functions}\label{ex:Bessel}
Let $\Psi(u)=u(u+\nu)$, $\nu {{ \geq}}  -\frac12$. We get that $
\WPn= n!\frac{\Gamma(n+\nu+1)}{\Gamma(\nu+1)}$
and thus
\begin{equation} \label{eq:bessel}
 \mathcal{J}_{\Psi}(t)=\sum^{\infty}_{n=0} (-1)^n\frac{\Gamma(\nu+1)}{\Gamma(n+\nu+1)} \frac{t^{2n}}{n!} = \Gamma(\nu+1) t^{-\nu} J_{\nu}\left( 2 t \right)
\end{equation}
where $J_{\nu}$ stands for  the  Bessel function of  order $\nu$. It is well-known that $\mathcal{J}_{\Psi} \in \De$.
\subsubsection{Confluent hypergeometric function} \label{ex:hyper}
Let $0<a<1<a+b$ and $\Psi(u)= u\frac{u+1-a}{u+b}$. Note that, in this case $\sigma=0$, and hence the support of the variable ${\rm{D}}_{\phi,\theta}$ is the real line. Moreover, simple algebra yields
\begin{eqnarray*}
\Psi(u) &=&u\frac{1-a}{b}+u\int^{\infty}_{0} \left(1-e^{-ur} \right) \left( a+b-1 \right) e^{-br} dr.
\end{eqnarray*}
We have $\WPn=n! \frac{\Gamma(b+1)\Gamma(n+2-a)}{\Gamma(n+b+1)\Gamma(2-a)}$
and thus
\begin{equation*}
 \mathcal{J}_{\Psi}(t)=\sum^{\infty}_{n=0} (-1)^n \frac{\Gamma(n+b+1)\Gamma(2-a)}{\Gamma(b+1)\Gamma(n+2-a)} \frac{t^{2n}}{n!} = {}_{1}F_{1}\left(b+1;2-a; -t^2\right)
\end{equation*}
where ${}_{1}F_{1}$ is the confluent hypergeometric function. If $a+b=2,3,\ldots$, then $\mathcal{J}_{\Psi} \in \De$, see e.g.~\cite[Theorem 4]{Hyp11}. Therein, the authors conjecture, in particular, that $\mathcal{J}_{\Psi} \in \De$ for all $a+b>1$.  %This result follows by an application of  Laguerre Theorem, see e.g.~\cite[Theorem 4]{Runck}.  %Note that, with $\phi_p(u)=\phi(u+p), p\in \mathbb N$,
%\begin{equation*}
 %\mathcal{I}_{\phi_p}(t)=  {}_{1}F_{1}\left(b+p+1;2+p-a; t\right).
%\end{equation*}

\subsubsection{Fox-Wright function}\label{sec:fw}
%\textcolor{red}{regarde dans mon article a refined factorization of the exponential law, l'exemple 1 pour les details de cet exemple j'ai juste fait une $T_{\beta-1/\alpha}$ transform de l'exemple}
Let $\alpha \in \left(0,1\right)$ and $\beta \geq \alpha$, then
\begin{equation*}
\Psi(u)=u\frac{\Gamma(\alpha u +\beta)}{\Gamma(\alpha (u-1) +\beta)} = u\frac{\Gamma(\beta)}{\Gamma(\beta-\alpha)}+\frac{u}{\Gamma(1-\alpha)} \int^{\infty}_{0}\left(1-e^{-ur} \right) \frac{e^{-\frac{\beta}{\alpha}r}}{\left(1-e^{-r/\alpha} \right)^{\alpha+1} }dr
\end{equation*}
yielding $
\WPn =n!  \frac{\Gamma(\alpha n +\beta)}{\Gamma(\beta)} $
and thus
\begin{equation*}
 \mathcal{J}_{\Psi}(t)=\sum^{\infty}_{n=0}(-1)^n \frac{\Gamma(\beta)}{\Gamma(\alpha n+\beta)} \frac{t^{2n}}{n!} = {}_0\Psi_1^{*} \left(  \left(\beta, \alpha\right); -t^2\right)
 \end{equation*}
 where ${}_0\Psi_1^{*}$ is the normalized Fox-Wright function, which is in $\De$ by Laguerre Theorem, see e.g.~\cite[Theorem 4]{Runck} as the mapping $z\mapsto \frac{1}{\Gamma(\alpha z+\beta)} \in \mathds{LP}_+$. % Moreover, with $\phi_p(u)=\phi(u+p), p\in \mathbb N$,
%\begin{equation*}
 %\mathcal{I}_{\phi_p}(t)=  {}_0\Psi_1^{*} \left(  \left(\beta + \alpha p, \alpha\right); t\right).
%\end{equation*}

\subsubsection{Mittag-Leffler function} \label{ex:ML}
Let now $\alpha \in \left(1,2\right)$ and $\beta \in \left(\alpha-1, \alpha\right)$ and set
\begin{equation*}
\Psi(u)=\frac{\Gamma(\alpha u +\beta)}{\Gamma(\alpha u +\beta-\alpha)}=\frac{\alpha \Gamma(\beta)}{\Gamma(\beta-\alpha+1)}u+ \frac{(\alpha-1)u}{\Gamma(2-\alpha)}\int^{\infty}_{0} \left(1-e^{-ur} \right) \frac{e^{-\frac{\beta}{\alpha}r}}{\left(1-e^{-r/\alpha} \right)^{\alpha}} dr.
\end{equation*}
Then, $
\WPn= \frac{\Gamma(\alpha n +\beta)}{\Gamma(\beta)} $ and
\begin{equation*}
 \mathcal{J}_{\Psi}(t)=\sum^{\infty}_{n=0} (-1)^n \frac{\Gamma(\beta)}{\Gamma(\alpha n+\beta)} t^{2n}= \Gamma(\beta) E_{\alpha,\beta}(-t^2)
 \end{equation*}
where $E_{\alpha,\beta}(t)=  \sum^{\infty}_{n=0} \frac{t^n}{\Gamma(\alpha n+\beta)}  $ is the general Mittag-Leffler function. It is well known, see \cite[Ch.~3.2]{zeroM}, that, for $ 1<\alpha <2$, $E_{\alpha,\beta}$ have nonreal zeros. %Moreover, with $\phi_p(u)=\phi(u+p), p\in \mathbb N$,
%\begin{equation*}
% \mathcal{I}_{\phi_p}(t)=\sum^{\infty}_{n=0}(-1)^n \frac{\Gamma(\alpha p +\beta)\Gamma( n+p+1)}{\Gamma( p+1)\Gamma(\alpha n+\alpha p + \beta)} \frac{t^{2n}}{n!} = {}_1\Psi_1^{*} \left( (p+1,1), \left(\beta+\alpha p, \alpha\right); -t^2\right)
% \end{equation*}
% where ${}_1\Psi_1^{*} \left( (p+1,1), \left(\beta+\alpha p, \alpha\right); t\right)$ is another Wright function.
For the last 3 examples, we do not know whether $\mathcal{J}_{\Psi} \in \De$ or not.
\subsubsection{The Barnes-Hypergeometric function}
Let now $\alpha \in \left(1,2\right)$,  $\rho \in (0,1/\alpha]$  and set
\begin{equation*}
\Psi(u)=\frac{\Gamma(\alpha \rho+\alpha u)u}{\Gamma(\alpha u) }= u \int^{\infty}_{0} (1-e^{-ur}) \frac{\rho}{\Gamma(1-\alpha \rho)} e^{-\rho r} \left( 1- e^{-r/\alpha} \right)^{-\alpha \rho-1} dr.%=\frac{\alpha \Gamma(\beta)}{\Gamma(\beta-\alpha+1)}+ \frac{\alpha-1}{\Gamma(2-\alpha)}\int^{\infty}_{0} \left(1-e^{-ur} \right) \frac{e^{-\frac{\beta}{\alpha}r}}{\left(1-e^{-r/\alpha} \right)^{\alpha}} dr.
\end{equation*}
Observe that
\begin{align}
\WPn&=n!\frac{G\left(1,\frac{1}{\alpha}\right)}{G\left(1+\rho,\frac{1}{\alpha}\right)}
\frac{G\left(n+1+\rho,\frac{1}{\alpha}\right)}{G\left(n+1,\frac{1}{\alpha}\right)} \label{eq:W_phm}
\end{align}
where the Barnes G-function $G(n;\frac{1}{\alpha})$ is defined as  the unique log-convex solution to recurrence equation $
G(n+1;\frac{1}{\alpha})=\Gamma\left(\alpha n\right)G(n;\frac{1}{\alpha})$, and refer to \cite{Loeffen} for more details on this example. Then,
\begin{equation*}
 \mathcal{J}_{\Psi}(t)=\sum^{\infty}_{n=0} (-1)^n \frac{G\left(n+1,\frac{1}{\alpha}\right)G\left(1+\rho,\frac{1}{\alpha}\right)}{G\left(1,\frac{1}{\alpha}\right)
 G\left(n+1+\rho,\frac{1}{\alpha}\right)}\frac{t^{2n}}{n!}.
 \end{equation*}
%Finally, with $\phi_p(u)=\phi(u+p), p\in \mathbb N$, one gets
%\begin{equation*}
% \mathcal{I}_{\phi_p}(t)=\sum^{\infty}_{n=0} (-1)^n \frac{G\left(n+1+p,\frac{1}{\alpha}\right)G\left(1+\rho+p,\frac{1}{\alpha}\right)}{G\left(1+p,\frac{1}{\alpha}\right)G\left(n+1+\rho+p,\frac{1}{\alpha}\right)}\frac{t^n}{n!}.
% \end{equation*}

\subsubsection{Hypergeometric function}
Let $\alpha>0$, and, writing  $\alpha_{1}= \frac{\alpha + \sqrt{\alpha^{2}+4}}{2}$ and $\alpha_{2}= \frac{\alpha - \sqrt{\alpha^{2}+4}}{2}$, consider
\begin{equation*}
\Psi(u)= \frac{u}{u+\alpha} \left( u + \alpha_{1} \right) \left( u + \alpha_{2} \right).
\end{equation*}
It follows that $
\WPn= n! \frac{\Gamma(\alpha+1)\Gamma(n+\alpha_{1}+1)\Gamma(n+\alpha_{2}+1)}{\Gamma(n+\alpha+1)\Gamma(\alpha_{1}+1)\Gamma(\alpha_{2}+1)} $ and thus
\begin{equation*}
 \mathcal{J}_{\Psi}(t)=\sum^{\infty}_{n=0} (-1)^n \frac{\Gamma(n+\alpha+1)\Gamma(\alpha_{1}+1)\Gamma(\alpha_{2}+1)}{\Gamma(\alpha+1)\Gamma(n+\alpha_{1}+1)\Gamma(n+\alpha_{2}+1)} \frac{t^{2n}}{n!} = {}_{1}F_{2}\left(\alpha+1;\alpha_{1}+1,\alpha_{2}+1; -t^2\right)
\end{equation*}
where ${}_{1}F_{2}$ is an hypergeometric function. %Moreover, with $\phi_p(u)=\phi(u+p), p\in \mathbb N$,
%\begin{equation*}
 %\mathcal{I}_{\phi_p}(t)= {}_{1}F_{2}\left(\alpha+p+1;\alpha_{1}+p+1,\alpha_{2}+p+1; -t^2\right).
%\end{equation*}
\subsubsection{Power-gamma function}
Let $\alpha \in (0,1)$, $\gamma\geq 0$ and consider
\begin{equation*}
\Psi(u)=  u\left( u + \gamma\right)^{\alpha} =\gamma^{\alpha}u+ u\int^{\infty}_{0} \left(1-e^{-ur} \right) e^{-\gamma r}\frac{r^{-\alpha-1}}{\Gamma(-\alpha)} dr.
\end{equation*}
It follows that $
\WPn= n!\frac{\Gamma^{\alpha}(n+\gamma+1)}{\Gamma^{\alpha}(\gamma+1)} $ and thus
\begin{equation*}
 \mathcal{J}_{\Psi}(t)=\sum^{\infty}_{n=0} (-1)^n\frac{\Gamma^{\alpha}(\gamma+1)}{\Gamma^{\alpha}(n+\gamma+1)}  \frac{t^{2n}}{n!}.
\end{equation*}
%More generally, $\Be$, the set of Bernstein functions, being  closed under composition, one gets, for instance, that, for any $\alpha \in (0,1)$ and  $\phi \in \Be$, $\phi^{\alpha} \in \Be$ as $\bar\phi(u)=  u ^{\alpha} \in \Be$ and one derives that
%\begin{equation*}
% \mathcal{I}_{\phi^{\alpha}}(t)=\sum^{\infty}_{n=0} (-1)^n\frac{1}{W_{\phi}^{\alpha}(n+1)}  \frac{t^{n}}{n!}.
%\end{equation*}

  \section{Proofs}\label{sec:proof}
  \subsection{Proofs of Section \ref{sec:Luk}}
  \subsubsection{Proof of Theorem \ref{it:dif}}
First, resorting to the Laguerre theory, see e.g.~Borel \cite[Ch.\ 2]{Borel},
on the additional number of zeros obtained by differentiating an entire function
of finite genus, we get  that if $\mathcal F\in \mathds{LP}$ then $t\mapsto\mathcal F^{(1)}(t)/t$ and $t\mapsto \mathcal F^{(2)}(t)$
are also in $\mathds{LP}$.
Recall that $\mathds{D}_L$ contains even and positive definite
functions in $\mathds{LP}$
that take value $1$ at $0$.
Suppose that $\mathcal F$ is even with $\mathcal F(0)=1$. it is easy to see that
$L^{(p)}\mathcal F$, $p=1,2$, are also even and take value $1$ at $0$.
Hence, to show that $L^{(p)}( \De)\subset \De$ it suffices to show
that if $\mathcal F$ is a characteristic function then so are the $L^{(p)}\mathcal F$, $p=1,2$.
One can check the latter by a probabilistic argument. First, recall
if ${\rm X}$ is a positive random variable with finite expectation $m_{\rm X}(1)$
and distribution $F_{\rm X}$ then the random variable ${\rm X}(1)$
with distribution
\[
F_{{\rm X}(1)}(dx) = \frac{x}{m_{\rm X}(1)} F_{\rm X}(dx), \: x\geq0,
\]
is called the size-biased version of ${\rm X}$. Size-biasing appears frequently and naturally
in probability theory, most notably in the theory of stationary point processes
on the real line and in the theory of branching processes and random walks.

%\tkside{You see, I avoid using $\mathbb P$\mbreak
%or $\mathbb E$ this way. I hope it's OK.}
\begin{lemma}\label{lem:Lm}
Let  $\mathcal F_X$ be the characteristic function of
a real-valued random variable ${\rm X}$
with finite second moment. Let ${{\rm X}^2(1)}$ be the size-biased version
of ${\rm X}^2$. Then $L^{(p)}\mathcal F_X$ is the characteristic function
of the random variable
\begin{equation}
\label{UX}
{\rm U}^{(p)} \times \sqrt{{\rm X}^2(1)}
\end{equation}
where ${\rm U}^{(p)}$ is a uniform
random variable on the interval $[-1,1]$ if $p=1$
or a random variable taking values $+1$ or $-1$ with probability $1/2$ each
if $p=2$, and, in both cases, is taken  independent of ${\rm X}^2(1)$.
\end{lemma}
\begin{proof}
We first recall  that, for a random variable $X$, we denote its distribution function by $F_X$. % is the distribution function of  ${\rm X}$, that is, $F(x)$ is the
Then, using the definition of size-biasedness and writing ${\rm X}^{(p)}={\rm U}^{(p)} \times \sqrt{{\rm X}^2(1)}$  we obtain, since the variables are independent, that
\[
F_{{\rm X}^{(1)}}(dx) = \frac{1}{m_{\rm X}(2)}\left(\int_{|y|>|x|} |y| F_{{\rm X}(1)}(dy) \right) dx,
\quad F_{{\rm X}^{(2)}}(dx) = \frac{x^2}{m_{\rm X}(2)} F_{{\rm X}(1)}(dx),
\]
where, by assumption,  $m_{\rm X}(2)
= \int_{-\infty}^\infty x^2 F_{{\rm X}}(dx)<\infty$.
We can then directly verify that
$L^{(p)}f(t) =  \int_{-\infty}^\infty e^{itx} F_{{\rm X}^{(p)}}(dx)$, $p=1,2$.
\end{proof}

\subsubsection{Proof of Theorem \ref{thm:DM}}
First, note that if $\mathcal{F}_{{\rm{X}}} \in \Pp_+$, then for any random variable
${{\rm{I}}}$, $\Lambda_{{\rm{I}}}\mathcal F_{{\rm{X}}} \in \Pp_+$ as it is the characteristic function of the random variable  ${\rm{X}}{\rm{I}}$ where the two random variables are considered to be independent. Next,  take  $\mathcal{F}_{{\rm{D}}} \in \Pp_+$ even and entire, that is, $\rm{D}$ is  a symmetric real-valued random variable such that, for any $t\in \R$,
  \begin{eqnarray*}
  \mathcal F_{{\rm{D}}}(t)&=&\sum^{\infty}_{n=0}(-1)^{n} \frac{\mathcal{M}_{{\rm{D}}}(2n) }{(2n)!}t^{2n}.
 \end{eqnarray*}
 Moreover, under the conditions of the item \eqref{it:Lag1}, that is $\varphi \in \mathds{LP}$ and $G \in \mathds{LP}_+$, Laguerre's theorem \cite[Theorem 4]{Runck} entails that  the function
\[
f(z)=\sum^{\infty}_{n=0}G(n)a_{\varphi}(n)z^{n}
\]
is an entire function with only real zeros.
We made the assumption that
$G(n)a_{\varphi}(n)={\mathcal{M}_{\rm{D}}(2n)\mathcal{M}_{\rm{I}}(2n)}/{(2n)!}$.
%and, ${\rm{D}}$ is independent of ${\rm{I}}$,
Since these are nonnegative numbers, $f$ cannot have nonnegative zeros.
On the other hand,
\begin{align*}
  \Lambda_{{\rm{I}}}\mathcal F_{{\rm{D}}}(t)= \int_{\R} \mathcal F_{{\rm{D}}}(xt)F_{{\rm{I}}}(dx)
&=\sum^{\infty}_{n=0}(-1)^{n} \frac{\mathcal{M}_{{\rm{D}}}(2n)\mathcal{M}_{{\rm{I}}}(2n) }{(2n)!}t^{2n}\\
  &=\sum^{\infty}_{n=0}(-1)^n G(n)a_{\varphi}(n)t^{2n}=f(-t^2)
\end{align*}
 where the interchange of the integral and sum is justified by a classical Fubini argument as the series defines an entire function.
Since $f$ has only real negative zeros, it follows that $\Lambda_{{\rm{I}}}\mathcal F_{{\rm{D}}}$ is an even entire function with only real zeros. Since $\Lambda_{{\rm{I}}}\mathcal F_{{\rm{D}}}$
is a characteristic function, it is in $\mathds P_+$. Hence
$\Lambda_{{\rm{I}}}\mathcal F_{{\rm{D}}} \in \mathds{D}_L$.
%\sout{This combines with Schoenberg Theorem, recalled in \eqref{Zeno}, show that $ \Lambda_{{\rm{I}}}\mathcal F_{{\rm{D}}} \in \D_L$ and hence  complete the proof of the first item.}

For the second one,  let $\mathcal{F} \in \De$, and thus, one has from \eqref{eq:pol}, that $\mathcal{F}$ is of the form \eqref{eq:defPL}, and thus $t\mapsto \mathcal{F}(\sqrt{t}) \in \mathds{LP}$.
Since $\mathcal{F} \in \De$
there exists a symmetric real-valued random variable $\rm{D}$
such that $\mathcal{F}=\mathcal{F}_{{\rm{D}}}$. Proceeding as above,   we get, for any $t\in \R$,
 \begin{eqnarray*}
  \Lambda_{{\rm{I}}}\mathcal F_{{\rm{D}}}(t)&=& \sum^{\infty}_{n=0}(-1)^{n} \frac{\mathcal{M}_{{\rm{D}}}(2n)\mathcal{M}_{{\rm{I}}}(2n) }{(2n)!}t^{2n}.
 \end{eqnarray*}
We apply again Laguerre's Theorem \cite[Theorem 4]{Runck}.
On the one hand, by assumption, there is a function $\varphi \in \mathds{LP}_+$
such that $\mathcal{M}_{{\rm{I}}}(2n)=\varphi(n)$.
On the other hand, we observed above that
\[
f(z)=\sum^{\infty}_{n=0}
\frac{\mathcal{M}_{{\rm{D}}}(2n)}
{(2n)!}z^{n}
= \mathcal F(\sqrt{z})  \in \mathds{LP}.
\]
Since $\Lambda_{{\rm{I}}}\mathcal F_{{\rm{D}}}(t) = f(-t^2)$ we conclude,
as above, that $\Lambda_{{\rm{I}}}\mathcal F_{{\rm{D}}} \in \mathds{D}_L$.

\subsubsection{Proof of Theorem \ref{thm:conv}}
Let first $(\mathcal F_{n})_{n\geq0}$ be a sequence in $\D$ and assume that for all $t$ in a bounded interval, $\lim_{n\to\infty}\mathcal F_{n}=\mathcal F$ uniformly. Then, according to the L\'evy continuity theorem, see \cite[Theorem 8.28]{Brei}, $\mathcal F \in \Pp_+$ and is continuous on $\R$. Let us write, for all $t \in \R, n\geq0$, $G_n(t)=\frac{1}{\mathcal F_n(t)}$. Since by assumption, for all $n\geq0$, $t\mapsto G_n(it)\in \Pp_+$, and, $\mathcal F_n$ is real, even and non-vanishing around $0$, we get that $G_n$ is well defined and even in a neighborhood of $0$. Hence, it is the moment generating function of a (unique) symmetric random variable. Moreover, we have, for all $t \in \R$,
$\lim_{n\to\infty}G_{n}(t)=G (t)=\frac{1}{\mathcal F(t)}$, and, by continuity of $\mathcal F $ and the fact that $\mathcal F(0)=1 $, there exits $\alpha>0$ such that $G$ is finite-valued on $|t|<\alpha$. According to \cite[Theorem 3]{curtis}, $G$ is the moment generating function of a random variable, whose law is uniquely determined by $G$, see \cite[Theorem 1]{curtis}. Finally, its characteristic function is plainly the mapping $t\mapsto G(it)$ and is such that, for all $t\in \R$, $\mathcal F(t) G(it)=1$, that is $\mathcal F \in \D$.

Now assume that the sequence $(\mathcal F_{n})_{n\geq0}$ is in $\De$. We shall show that its locally uniform limit  $\mathcal F\in  \De$.
Since, for all $n \geq0$,  $\mathcal F_n\in  \De$, by  Theorem \ref{thm:sch}, we have that $\mathcal G_n: t\mapsto\frac{1}{\mathcal F_n(it)}$ is the characteristic function of a symmetric P\'olya frequency function. As $\De \subset \D$, from the previous proof, we deduce that, for all $t \in \R$, $\lim_{n\to\infty}\mathcal G_n(t)=\mathcal G(t)=\frac{1}{\mathcal F(it)}$ with $\mathcal G \in \Pp_+$. However, the set of P\'olya frequency functions (probability measures)  is closed under weak convergence as it is  the closure, under weak convergence, of probability measures whose characteristic functions are the reciprocal of a polynomial having only real roots, see \cite[Chap.~III, Theorem 4.1]{Wid}. We obtain, using  the fact that $\mathcal G_{n}$ is real and even for all $n\geq0$,  that the mapping $ t\mapsto \mathcal{F}(t)=\frac{1}{\mathcal G (it)}\in \LPnew$  and hence  ${\mathcal F} \in \Pp_+ \cap \LPnew=\De$ which completes the proof.
\subsection{Proofs of Section \ref{intro part 2}}
\subsubsection{A Wiener-Hopf type mapping between sets of negative-definite functions}
  We start the proof by   a one-to-one mapping, emanating from the Wiener-Hopf factorization,  between the set $\N_{\D}$ and a set of Bernstein functions.
This allows us to  provide a representation of the  coefficients of the function
in terms of a Stieltjes moment sequence that will be helpful in identifying
one of the random variables solving the van Dantzig problem.

To this end we introduce
some
notation.
First, let us denote by $\Be$  the set of Bernstein functions,
i.e.\ functions $\phi: (0,\infty) \to [0,\infty)$ having derivatives of all orders such that
$(-1)^{n+1} \phi^{(n)} \ge 0$ for all $n \ge 1$;
see \cite[Ch.\ 3]{SchillingSongVondracek10}
These functions are in a one-to-one correspondence with functions that
admit the so-called L\'evy-Khintchine representation
\begin{equation} \label{eq:bern}
\phi(u) =  \nu+ \frac{\sigma^2}{2}u +\int_{0}^{\infty} (1-e^{-ur}) \mu(dr),\ u\geq 0,
\end{equation}
where $\nu,\sigma\geq0$ and  $\mu \in \mathbb{M}^+(\R^+)$ such that $\int^{\infty}_{0} \inf(1, r) \mu(dr) <\infty$.
The set $\Be$ is invariant under several transformations.
In particular, it is a convex cone, and, it is also stable by the action of the semigroup of translations, Indeed, for any $\beta \geq 0$, easy algebra yields
\begin{equation} \label{eq:bernz}
 \phi(u+\beta) = \phi(\beta) +  \frac{\sigma^2}{2}u +\int_{0}^{\infty} (1-e^{-ur})e^{-\beta r} \mu(dr), \  u\geq 0,
 \end{equation}
 which is plainly a Bernstein function as $\phi(\beta)\geq 0$ and $e^{-\beta r} \mu(dr)$ is a L\'evy measure as defined above.
%$\phi^{(1)}$ is completely monotone and by
 We shall also need the subset $\Be_J \subset \Be$ which is the convex cone  of Bernstein functions which take the form
   \begin{equation} \label{eq:bernJ}
 \phi(u) =   \nu+ \frac{\sigma^2}{2}u +\int_{0}^{\infty} (1-e^{-ur}) \overline{\mu}(r)dr
 \end{equation}
 where ${\nu,\sigma \geq 0}$ and  $\overline{\mu}$ is a non-negative and non-increasing function on $\R^+$ such that $\int^{\infty}_{0} \inf(1,r)\overline{\mu}(r)dr <\infty$. We refer to the  monograph \cite{SchillingSongVondracek10} for an excellent account on all these sets of functions.

Next, for a Radon measure $\mu$ on $(0,\infty)$ and $\theta\ge 0$, let
\[
\overline\mu_\theta(r):= \int_r^\infty e^{\theta(r-s)} \mu(ds), \quad r >0.
\]
Letting  $\mathcal M_\theta$ be the set of all functions $\overline\mu_\theta$,
where $\mu$ ranges over Radon measures, such that
$\int_0^\infty \min(1,r) \overline \mu_\theta(r) dr<\infty$,
we define, in relation to the van Dantzig problem, the following set
of Bernstein functions:
 \begin{equation}\label{eq:defBD}
  \Be_{\D}=\left\{\phi \in \Be; \mu(dr) =\overline{\mu}_{\theta}(r)dr, \textrm{ with } \overline{\mu}_{\theta}(r) \in  \mathcal{M}_{\theta}, 0\leq \theta {{ \leq}}  \frac12 \right\}.
\end{equation}
%where  \begin{equation*}\label{eq:defM}
%         \mathcal{M}_{\theta}=\left\{\overline{\mu}_{\theta}; \overline{\mu}_{\theta}(r)= \int_{r}^{\infty}e^{-\theta (s-r)}\mu(ds),  \mu \in \mathbb{M}^+(\R^+) \textrm{ and }  \int^{\infty}_{0} \inf(1, r) \overline{\mu}_{\theta}(r)dr <\infty \right\}.
%       \end{equation*}
Note first that when $\theta>0$ the function $\overline \mu_\theta$ may
fail to be monotone, for a suitable choice of $\mu$.
Note also that when $\theta =0$ we have
 $\overline{\mu}_{0}(r)= \int_{r}^{\infty}\mu(ds)$ is simply a
non-negative and non-increasing function that satisfies the
integrability condition above,meaning that $\Be_J$ is a strict subset of $\Be_{\D}$.
%, meaning that  $\Be_J \subsetneq \Be_{\D}$.

 We have the following.
       \begin{proposition}\label{prop:WH}
       \begin{enumerate}
       \item  There exists a one-to-one mapping between the sets $\N_{\D}$ and $\Be_{\D}$. More specifically,  for any  $ \Psi \in \N_{\D}$ of the form \eqref{LLK}, we have, for any $u\geq0$,
       \begin{equation}\label{eq:WH}
       \Psi(u)=(u-\theta)\phi(u)
       \end{equation}
        where, with $\nu_{\theta}=\frac{\kappa}{\theta}\mathbb{I}_{\{\theta>0\}}+\Psi^{(1)}(0^+)\mathbb{I}_{\{\theta=0\}}$, we have set $\phi(u)=\nu_{\theta}+\frac{1}{2}\sigma^{2}u+\int_0^{\infty}(1-\mathrm{e}^{-u r})\overline{\mu}_{\theta}(r)dr \in \Be_{\D}$. Moreover, with such a notation, we have, for all $z\in \C$,
 \begin{equation}\label{eq:defIpz}
\mathcal{J}_{\Psi}(z)=\mathcal{J}_{\phi,{\theta}}(z)  =\sum^{\infty}_{n=0} \frac{(-1)^n\Gamma(1-\theta)}{\Gamma(n+1-\theta)W_{\phi}(n+1)}z^{2n}.
  \end{equation}
  \item \label{it:tb}Let us define the mapping $\mathcal{T}_\beta f(u)=\frac{u}{u+\beta}f(u+\beta)$. Then, $(\mathcal{T}_\beta)_{\beta\geq0}$ form a semigroup,
 i.e., $\mathcal{T}_\beta \circ \mathcal{T}_\gamma
= \mathcal{T}_{\beta+\gamma}$,
and, for any  $\beta\geq 0$, $\mathcal{T}_\beta(\N)\subset \N$ and  $\mathcal{T}_\beta(\N_\D)\subset \N_{\D}$.
       \item\label{it:tb2} Let us write $\overline{\mathcal{T}}_\beta f(u)=\frac{u-\beta}{u+\beta}f(u+\beta)$. Then, for any  $\beta\geq \theta$, $\overline{\mathcal{T}}_\beta(\N)\subset \N$ and for any $\theta\leq \beta {{ \leq}}  \frac12$, $\overline{\mathcal{T}}_\beta(\N_\D)\subset \N_\D$. %(resp.~ $u\mapsto \widetilde{\mathcal{T}}\Psi(u)=\frac{u-\frac{1}{2}}{u+\frac{1}{2}}\frac{u}{u+1}\Psi(u+1) \in \N_{\D}$)
\end{enumerate}
       \end{proposition}
  \begin{proof}
 Let $ \Psi \in \N_{\D}$.  Then,  by means of the Wiener-Hopf factorization of L\'evy-Khintchine exponents of spectrally negative L\'evy process, see e.g.~\cite{Kyprianou06}, there exists $\phi \in \Be$ such that, for all $u\geq0$,
\begin{equation}\label{eq:WH}
  \Psi(u)=(u-\theta)\phi(u)
\end{equation}
where we recall that $\theta=\sup\{u\geq0; \Psi(u)=0\}$, see \eqref{eq:zeta}.
In order to characterize $\phi$, we write $\overline{\nu}_{\theta}=a+\sigma^{2}\theta+\int_0^{1}(1-\mathrm{e}^{-\theta r})r\mu ({\rm d}r)-\int_1^{\infty} e^{-\theta r}\overline{\mu}_{\theta}(r)dr$ and,  observe first that
\begin{eqnarray*}%\label{eq:WH}
  \Psi(u+\theta)&=& (\overline{\nu}_{\theta}+\int_1^{\infty} e^{-\theta r}\overline{\mu}_{\theta}(r)dr)u +\frac{1}{2}\sigma^{2}u^{2}-\int_0^{\infty}(1-\mathrm{e}^{-u r}-u r\mathbb{I}_{\{r<1\}})e^{-\theta r}\mu ({\rm d}r)\\
  &=&  u\left(\overline{\nu}_{\theta}+\frac{1}{2}\sigma^{2}u+\int_0^{\infty}(1-\mathrm{e}^{-u r}) e^{-\theta r}\overline{\mu}_{\theta}(r)dr\right)
\end{eqnarray*}
and, thus
\begin{eqnarray*}%\label{eq:WH}
  \Psi(u)&=&  (u-\theta)\left(\nu_{\theta}+\frac{1}{2}\sigma^{2}u+\int_0^{\infty}(1-\mathrm{e}^{-u r}) \overline{\mu}_{\theta}(r)dr\right)\\
\end{eqnarray*}
which, after recalling that $-\theta\phi(0)=\Psi(0)=-\kappa$ and $\Psi^{(1)}(0^+)=\phi(0)$ if $\theta=0$, completes the proof of the first item as the rest follows at once. Next,  the claim of item \eqref{it:tb} in the case $\beta=0$ is obvious and thus we assume that $\beta>0$. In \cite{CKP}, see also \cite[Proposition 2.1]{Patie-Savov2}, it is shown that $(\mathcal{T}_\beta)_{\beta\geq0}$ is a semigroup and the set $\N$ is invariant under the action of the mapping $\mathcal{T}_\beta, \beta\geq0$, that is $\mathcal{T}_\beta(\N)\subset \N$. Thus, it simply remains to show that $\mathcal{T}_\beta(\Nd)\subset \Nd$ for all $\beta\geq0$. Let $\Psi \in \Nd$. Then, note that $\mathcal{T}_\beta \Psi(0)=0$, i.e.~$\kappa=0$ in \eqref{LLK}, and, there exists $0\leq \theta {{ \leq}}  \frac12$ such that $\Psi(\theta)=0$. Suppose,  first, that $\beta\geq \theta$, then recalling that $\Psi$ is non-decreasing on $[\theta,\infty)$ with $\Psi(\theta)\geq0$, we get the claim, from \eqref{eq:propNd}, by observing that $(\mathcal{T}_\beta \Psi)^{(1)}(0)=\frac{\Psi(\beta)}{\beta}\geq 0$. Next, assuming that $0<\beta <\theta$, one easily sees that $\mathcal{T}_\beta \Psi(\theta-\beta)=\frac{\theta-\beta}{\theta}\Psi(\theta)=0$ and since $0<\theta-\beta<\frac12$, we deduce  that also $\mathcal{T}_\beta \Psi \in \N_\D$, completing the proof of this item as the semigroup property is obvious. For the last claim, first  observe that,  $\overline{\mathcal{T}}_\beta \Psi(u)=\frac{u+\beta-2\beta}{u+\beta}\Psi(u+\beta)+\Psi(\beta)-\Psi(\beta)=\mathcal{T}_{2\beta,\beta} \Psi(u)-\Psi(\beta)$ (the last identity serves to fix a notation), and then we know from \cite[Proposition 2.2]{CKP}, that $\mathcal{T}_{2\beta,\beta}(\N)\subset \N$ and we obtain the statement since $\Psi(\beta)\geq 0$ as we choose $\beta\geq \theta$. Finally, if $\Psi \in \N_\D$, then, from \eqref{eq:propPsi}, $\Psi$ is {{ non-negative }}  on $[\frac{1}{2},\infty)$ and, from the preceding discussion, $\overline{\mathcal{T}}_{\beta} \Psi \in \N$ with $\overline{\mathcal{T}}_{\beta} \Psi(\frac{1}{2})=\frac{\frac{1}{2}-\beta}{\frac{1}{2}+\beta}\Psi(\frac{1}{2}+\beta){{ \geq  }}0$.
\end{proof}
%Finally,
%where it is not difficult to check, using \eqref{eq:defP}, that $\phi \in \Be$, and, when  $\theta=0$ then $\phi \in \Be_J$
%in relation to the van Dantzig problem,  we introduce the following subset of Bernstein functions
%\begin{equation}\label{eq:defBD}
%  \Be_{\D}=\{\phi \in \Be;\: \phi(u)=\frac{\Psi(u)}{u-\theta}, \Psi \in \N_\D\}.
%\end{equation}

\subsubsection{Some analytical properties of the entire functions in the class $\Dp$} \label{sec:PropJ}
\begin{lemma} \label{lem:L1}
Let $\Psi \in \Nd$, and, for $p=1,2$, we set  $L_1^{(p)}=L^{(p)}$, where we recall that $L^{(p)}$ was defined in \eqref{eq:lukacs},  and, for  $k=1,2,\ldots$, we define
\[ L_{k+1}^{(p)}=L^{(p)} \circ L_k^{(p)}.\]
\begin{enumerate}
  \item For any $k=1,2,\ldots$, we have
\begin{equation}\label{eq:propT}
    L_k^{(1)}\mathcal{J}_{\Psi}= \mathcal{J}_{\mathcal{T}_k\Psi}
\end{equation}
where we recall that the mapping $\mathcal{T}_k$ was defined in Proposition \ref{prop:WH}\eqref{it:tb}.
\item Moreover, for any $k=1,2,\ldots$, we have
\begin{equation}\label{eq:propT2}
   L_k^{(2)}\mathcal{J}_{\Psi}=  \mathcal{J}_{\overline{\mathcal{T}}^{2k}_{\frac{1}{2}}\Psi}
\end{equation}
where, $\overline{\mathcal{T}}^1_{\beta}=\overline{\mathcal{T}}_{\beta}$ and $\overline{\mathcal{T}}^k_{\beta} =\overline{\mathcal{T}}_{\beta}\circ \overline{\mathcal{T}}^{k-1}_{\beta}$, and,  the mapping $\overline{\mathcal{T}}_{\beta}$ was defined in Proposition \ref{prop:WH}\eqref{it:tb2}. %\Psi(u)=\frac{u-\beta}{u+\beta}\frac{u}{u+2\beta}\Psi(u+2\beta)$.
\end{enumerate}
\end{lemma}
\begin{proof}
 To prove the identity \eqref{eq:propT}, we first let  $k=1$ and observe that
\begin{equation}\label{eq:difJ}
    \mathcal{J}^{(1)}_{\Psi}(t)=\sum^{\infty}_{n=1} \frac{2n(-1)^n t^{2n-1} }{W_{\Psi}(n+1)}=\frac{-2t}{\Psi(1)}\sum^{\infty}_{n=0} \frac{(-1)^n t^{2n} }{\prod_{k=1}^{n}\frac{k}{k+1}\Psi(k+1)}= \frac{-2t}{W_{\Psi}(2)} \mathcal{J}_{\mathcal{T}_{1}\Psi}(t).
\end{equation}
Moreover, differentiating the right-hand side one more time yields
\[
\mathcal{J}^{(2)}_{\Psi}(0)=\lim_{t\to 0}\frac{-2}{W_{\Psi}(2)}\mathcal{J}_{\mathcal{T}_1\Psi}(t)=\frac{-2}{W_{\Psi}(2)}
\]
where we used that $\mathcal{J}^{(1)}_{\mathcal{T}_1\Psi}(0)=0$,  which, itself, follows by applying \eqref{eq:difJ} to $ \mathcal{J}_{\mathcal{T}_1\Psi}$. Hence $L^{(1)}\mathcal{J}_{\Psi}= \mathcal{J}_{\mathcal{T}_1\Psi}$. To complete the proof for $p=1$, we resort to an induction argument combined with the semigroup property of the mapping $\mathcal{T}_k$. Indeed, for any $k=1,\ldots$, one has
\begin{equation*}%\label{eq:propT}
    L_{k+1}^{(1)}\mathcal{J}_{\Psi}= L^{(1)}\circ L_{k}^{(1)}\mathcal{J}_{\Psi}= L^{(1)}\mathcal{J}_{\mathcal{T}_k\Psi}=\mathcal{J}_{\mathcal{T}_{k+1}\Psi}.
\end{equation*}
On the other hand, we note that
\begin{align*}%\label{eq:diffJ}
    \mathcal{J}^{(2)}_{\Psi}(t)
&=\frac{-2}{W_{\Psi}(2)}\frac{d}{dt}\sum^{\infty}_{n=0} \frac{(-1)^n  }{W_{\mathcal{T}_{1}\Psi}(n+1)}t^{2n+1}
    =\frac{-2}{W_{\Psi}(2)}\sum^{\infty}_{n=0} \frac{(-1)^n (2n+1)  }{W_{\mathcal{T}_{1}\Psi}(n+1)}t^{2n}
\\
&=\frac{-2}{W_{\Psi}(2)}\sum^{\infty}_{n=0} \frac{(-1)^n  }{\prod_{k=1}^{n}\frac{k-\frac{1}{2}}{k+\frac{1}{2}}\frac{k}{k+1}\Psi(k+1)}t^{2n}
\\
&=\frac{-2}{W_{\Psi}(2)} \mathcal{J}_{\overline{\mathcal{T}}^2_{\frac{1}{2}}\Psi}(t)
\end{align*}
which provides the proof of the last claim for $k=1$. As above invoking an induction argument yields
\begin{equation*}%\label{eq:propT}
    L_{k+1}^{(2)}\mathcal{J}_{\Psi}= L^{(2)}\circ L_{k}^{(2)}\mathcal{J}_{\Psi}= L^{(2)}\mathcal{J}_{\overline{\mathcal{T}}^{2k}_{\frac{1}{2}}\Psi}=\mathcal{J}_{\overline{\mathcal{T}}^{2k+2}_{\frac{1}{2}}\Psi}
\end{equation*}
which completes the proof.
\end{proof}

We proceed with the following result regarding the order and the type of the entire function $\mathcal{J}_{\Psi}$. It is a slight refinement  of a result obtained recently by Bartholm\'e and Patie \cite{Bart-Patie} in the case $\theta=0$, and, for the entire function $\mathcal{I}_{\Psi}(\sqrt{t})=\mathcal{J}_{\Psi}(i\sqrt{t})$, this later transformation affects simply the order by a factor of $2$. To deal with the general case $\theta>0$,  it is not difficult to see that both the order and type  remain the same when replacing the term $n!$ by $C\Gamma(n+1-\theta), C>0,$ in the coefficients of a power series. We left the details of the easy modification to the reader and recall that a measurable function $\ell$ on $\R^+$ is said to be slowly varying at infinity if, for every $u>0$, $\lim_{t\to\infty}\frac{\ell(ut)}{\ell(t)}=1$.
\begin{proposition}\cite[Proposition 2.1]{Bart-Patie} \label{prop:order}
Let $\Psi \in \Nd$. Then, $\mathcal{J}_{\Psi}$ is an entire function of order
\begin{equation}\label{eq: order of I}
\rho_{\Psi}=\frac{2}{\underline{\Psi}} \in \left[1,2\right]
\end{equation}
where $
\underline{\Psi} = \sup \{a>0; \: \lim_{u \to \infty} u^{-a}\Psi(u) = \infty \} =  \liminf_{u \rightarrow \infty} \frac{\ln \Psi(u) }{ \ln u}\in [1,2]$, is its so-called Blumenthal-Getoor lower index. Moreover, its type $\tau_{\Psi}$ is given by \begin{equation}\label{eq: type of I}
\tau_{\Psi}=\underline{\Psi}e^{-\frac{\underline{\Psi}-1}{\underline{\Psi}}}\limsup_{n \to \infty}e^{ -\frac{1}{\underline{\Psi}}\left(\frac{\int_{1}^{n}\ln \phi(u)du}{n}-(\underline{\Psi}-1)\ln n \right)} \geq \left(\limsup_{n \to \infty}\frac{n^{\underline{\Psi}-1}}{\phi(n)}\right)^{\frac{1}{\underline{\Psi}}}
\end{equation}
where we recall that $\Psi(u)=(u-\theta)\phi(u)$.
 In particular, $\rho_{\Psi}=1$ if $\Psi(u)=\frac{u^{2}}{\ell(u)}$, and, otherwise, $\rho_{\Psi}=\frac{2}{\mathfrak{a}}, \mathfrak{a} \in \left[ 1,2\right),$ if  $\Psi(u)=u^{\mathfrak{a}}\ell(u)$  where, in both cases, $\ell$ is a slowly varying function, and,
 $\rho_{\Psi}=2$ also if  $\int_{0}^{\infty}r\mu(dr)<\infty$.
\end{proposition}

\subsubsection{Useful bounds for the Bernstein-gamma functions}
 For a Bernstein function $\phi \in \Be$, we write $z_{\phi}=\inf\{u>0; \phi(-u)=0\} \in [0,\infty]$ and  $e_{\phi}=\sup\{u>0; \int_{1}^{\infty}e^{ur}\mu(dr)<\infty\} \in [0,\infty]$ and set
\[ a_{\phi}=\min(z_{\phi},e_{\phi}). \]
We also denote by    $W_{\phi}$ the so-called Bernstein-gamma function which admits the following  generalized Weierstrass product representation \begin{equation*}
W_{\phi}(z)=\frac{e^{-\gamma_{\phi}z}}{\phi(z)} \prod^{\infty}_{k=1} \frac{\phi(k)}{\phi(k+z)} e^{\frac{\phi'(k)}{\phi(k)}z}
\end{equation*}
where
\begin{equation*}
\gamma_{\phi} = \lim_{n \rightarrow \infty} \left( \sum^{n}_{k=1} \frac{\phi'(k)}{\phi(k)} - \log \phi(n) \right).
\end{equation*}
From \cite[Theorem 4.1]{PS-B}, we know that $W_{\phi}$ defines a zero-free and analytic function on the right-half plane $\Re(z)>-a_{\phi}$, which admits a meromorphic extension to the right-half plane $\Re(z)>-e_{\phi}$.
We mention that when $\phi(z)=z$, then $W_{\phi}$ boils down to the gamma function as the infinite product above corresponds to its Weierstrass representation and $\gamma_{\phi}$ is the Euler-Mascheroni constant. This class of functions has been thoroughly studied in the papers by Patie and Savov \cite{Patie-Savov, PS-B}. Thereout,  it has been shown that $W_{\phi} \in \Pp_+$ and it is the unique element in $\Pp_+$ solution to the functional equation $ W_{\phi}(z+1)={\phi}(z)W_{\phi}(z), W_{\phi}(1)=1$. We  also find  in \cite{PS-B} the following result which will be  useful in the sequel.
\begin{proposition}\cite[Proposition 6.2]{PS-B} \label{prop:condClass1}
Let $\phi\in\Be$ and write $\overline{\mu}(r)=\int_{r}^{\infty}\mu(ds), r>0$. Then, the following estimates holds.	
\begin{enumerate}
	  \item If $\sigma^2>0$  %then \eqref{eq:condClass1} holds for any $u< \frac{2}{\sigma^2}\left(\phi(0)+\overlin{\mu}(0)\right) \in (0,\infty]$. In fact, if $\mubr\!\lbrb{0}<\infty$ the limit in \eqref{eq:condClass1} is infinity for all $u>\frac{1}{\dr}\lbrb{\phi(0)+\mubr\!\lbrb{0}}$. Finally, regardless of the value of $\mubr(0)$,
we have, for any  $\epsilon, a>0$ such that $\int_{0}^{\infty} e^{-ar}\overline{\mu}(r)dr<1$,  as $|b|\to\infty$,
	\begin{equation}\label{eq:condClass11}
	\left|\frac{\Gamma(a+ib)}{W_{\phi}(a+ib)}\right| \leq C |b|^{-\frac{2}{\sigma^2}\left(\overline{\mu}(1/|b|)+\phi(0)\right)+\epsilon}
	\end{equation}
where $C>0$.
	  \item If $\sigma^2=0$ then, for any $ u\geq 0$ and $a>0$ fixed,
	\begin{equation}\label{eq:condClass1}
	\lim_{|b|\to \infty}|b|^u\left|\frac{\Gamma(a+ib)}{W_{\phi}(a+ib)}\right|=0.
	\end{equation}
	\end{enumerate}
	
\end{proposition}

\subsubsection{Characterization and properties of the van Dantzig pair of random variables}		\label{charprop}
 To define the van Dantzig pair of random variables that appear in Theorem \ref{thm:main}, we need to introduce some notations. Let $\Psi \in \N$ and recall that it is the Laplace exponent of a spectrally negative L\'evy process ${\rm{Y}}$, see \eqref{eq:LeY}. Then,
 according to Lamperti \cite{Lamperti}, there exists a  $1$-self-similar Markov process on $(0,\infty)$ denoted by ${\rm{X}}=({\rm{X}}_t)_{t\geq 0}$ such that
 \begin{equation}\label{eq:lamp}
   {\rm{X}}_t = e^{\mathrm{Y_{A_t}}}, \quad 0\leq t < \xi=\inf\{t>0;\: {\rm{X}}_t=0\},
 \end{equation}
 where $\mathrm{A}_t=\inf\left\{ s>0; \: \int_{0}^{s}e^{\mathrm{Y}_u}du>t\right\}$. Moreover,  if $\theta=0$, i.e.~the conditions \eqref{eq:zeta} are not fulfilled,  then $\xi=\infty$ almost surely and ${\rm{X}}$ is in fact a Feller process on $[0,\infty)$. On the other hand, if $0<\theta<1$, then $\xi$ is finite almost surely but there exists an unique  $1$-self-similar extension of the process ${\rm{X}}$ which is also a Feller process on  $[0,\infty)$ and  with the path property to leave the recurrent boundary point $0$ continuously, see \cite{Rivero}.
 From now on, if $0\leq \theta<1$, we denote by   $\mathrm{X}(\Psi)=(\mathrm{X}_t(\Psi))_{t\geq 0}$ the realization of the   $1$-self-similar Feller semigroup $(P^{\Psi}_t)_{t\geq0}$  on $[0,\infty)$ (that is the recurrent extension if $\theta>0$) associated to $\mathrm{Y}$  and
 let ${\rm{T}}_{\Psi}$ be a random variable whose distribution is that of
%\[ {\rm{T}}_{\Psi}=\inf\{ t>0;\: X_0(\Psi)=0 \textrm{ and }   X_t(\Psi)\geq 1\}\]
\[
\inf\{ t>0;\: \mathrm{X}_t(\Psi)\geq 1\},
\]
conditional on $\mathrm{X}(\Psi)$ starting from $0$.
%for its first passage time to the level $1$ when the process starts from $0$.
We also need to introduce the  exponential functional of a subordinator, a positive random variable which has been intensively studied over the last two decades and refer to \cite{Patie-Savov2} for a nice historical account  and a thorough study. Let $\phi \in \Be$ and let
\begin{equation}\label{eq:moX}
{\rm{I}}_\phi=\int_{0}^{\infty}e^{-{\rm{Z}}_t}dt
\end{equation}
where $({\rm{Z}}_t)_{t\geq 0}$ is a subordinator such that $ \mathcal{F}_{{\rm{Z}}_t}(iu)=e^{-\phi(u)t}, u,t\geq0$. Next, we recall, from \cite{CPY}, the following expression for the integer moments of ${\rm{I}}_{\phi}$,
\begin{equation}\label{eq:moX}
\mathcal{M}_{{\rm{I}}_{\phi}}(n)=\frac{n!}{W_{\phi}(n+1)}, \quad n= 0,1,\ldots,
\end{equation}
and the notation of the associated Markov operator, see  \eqref{eq:defLI},
\begin{equation}\label{eq:MI}
\Lambda_{{\rm{I}}_{\phi}} f(t)=\int_{0}^{\infty}f(xt)F_{{\rm{I}}_{\phi}}(dx).
\end{equation}
We are now ready to state the following result which can be found  in \cite{Patie}.
To emphasize the role played by the concept of intertwining in the van Dantzig problem,  we  provide  another original  proof.
Intertwining, in our context, refers to the relation \eqref{intertw} below.
To state the theorem, recall the notation
\[
\mathcal J_{-\theta}(x) = \Gamma(1-\theta) t^{\theta} J_{-\theta}(2x), \quad
\mathcal I_{-\theta}(x) = \Gamma(1-\theta) t^{\theta} I_{-\theta}(2x),
\]
appearing in \eqref{eq:bessel},
where $J_{-\theta}$ (respectively $I_{-\theta}$)  is the ordinary (respectively modified) Bessel function of the first kind of order $-\theta$.
%Recall that $\mathcal J_{-\theta}(x)$ appeared in \eqref{eq:bessel} where it was observed that, when $\theta < 1/2$,
%the pair $\mathcal J_{-\theta}(x), 1/\mathcal I_{-\theta}(x)$ is a van Dantzig pair with both elements in $\mathds D_L$.

\begin{theorem} \cite[Theorem 2.1]{Patie}\label{thm:patie}
 Let  $\Psi \in \N$ with $0\leq \theta<1$.
 Then, writing $\widetilde{\mathcal{I}}_{\Psi}(x)=\mathcal{I}_{\Psi}(\sqrt{x})$
and $\widetilde{\mathcal I}_{-\theta}(x) = \mathcal I_{-\theta}(\sqrt{x})$,
we have
   \begin{equation}\label{eq:IntI}
   \Lambda_{{\rm{I}}_{\phi}}\widetilde{\mathcal{I}}_{-\theta}(x)=\widetilde{\mathcal{I}}_{\Psi}({x}), \quad x>0,
   \end{equation}
 where  $\phi(u)=\frac{\Psi(u)}{u-\theta} \in \Be$.
%and $\widetilde{\mathcal{I}}_{-\theta}(x)={\mathcal{J}}_{-\theta}(i\sqrt{x})$,
%see \eqref{eq:bessel},  is the modified Bessel function of order $-\theta$.
Moreover,  ${\rm{T}}_{\Psi}$ is a positive self-decomposable random variable, and, for any $u>0$, we have
 \begin{eqnarray}\label{eq:LTT}
\mathcal F_{{{\rm{T}}_{\Psi}}}(iu)&=&\frac{1}{\mathcal{I}_{\Psi}(\sqrt{u})}=e^{-\phi_{\Psi}(u)}
 \end{eqnarray}
 where $\phi_{\Psi}$ is the Bernstein function defined in Theorem \ref{thm:main}.
\end{theorem}
\begin{proof}
 Let  $\Psi \in \N$ with $0\leq \theta<1$. We recall that in \cite{Patie-Savov} for the case $\theta=0$, and in \cite{PSZ}  for the case $0<\theta<1$, the following intertwining relation has been identified
 \begin{equation} \label{intertw}
P^{\Psi}_t\Lambda_{{\rm{I}}_{\phi}} = \Lambda_{{\rm{I}}_{\phi}}P^{\Psi_{\theta}}_t
 \end{equation}
 where  $\Psi_{\theta}(u)=u(u-\theta)\in \N$.
Then, using successively   Tonelli theorem and the identity \eqref{eq:moX},
we obtain
that, for any $x>0$,
\begin{eqnarray*}
\Lambda_{{\rm{I}}_{\phi}} \widetilde{\mathcal{I}}_{-{\theta}}(x)&=&\sum_{n=0}^{\infty}\frac{\mathcal{M}_{{\rm{I}}_{\phi}}(n)\Gamma(1-\theta)}{n!\Gamma(n+1-\theta)}x^n
=\sum_{n=0}^{\infty}\frac{\Gamma(1-\theta)}{W_{\phi}(n+1)\Gamma(n+1-\theta)}x^n=\widetilde{\mathcal{I}}_{\Psi}({x}).
 \end{eqnarray*}
 It is well-known that the mapping $x\mapsto d_q\widetilde{\mathcal{I}}_{-{\theta}}(x)=\widetilde{\mathcal{I}}_{-{\theta}}(qx), q>0,$ is a $q$-invariant for the squared Bessel semigroup $P^{\Psi_{\theta}}$, that is, for all $t\geq0$, $e^{-qt}P^{\Psi_{\theta}}_td_q\widetilde{\mathcal{I}}_{-{\theta}}(x)=d_q\widetilde{\mathcal{I}}_{-{\theta}}(x)$. Now observing that $d_q\Lambda_{{\rm{I}}_{\phi}}=\Lambda_{{\rm{I}}_{\phi}}d_q$, we deduce from the intertwining relation above that
  \begin{eqnarray*}
P^{\Psi}_td_q\widetilde{\mathcal{I}}_{\Psi}(x)&=&P^{\Psi}_td_q\Lambda_{{\rm{I}}_{\phi}}\widetilde{\mathcal{I}}_{-{\theta}}(x)= P^{\Psi}_t\Lambda_{{\rm{I}}_{\phi}}d_q\widetilde{\mathcal{I}}_{-{\theta}}(x) \\&=& \Lambda_{{\rm{I}}_{\phi}}P^{\Psi_{\theta}}_td_q\widetilde{\mathcal{I}}_{-{\theta}}(x)=e^{qt}\Lambda_{{\rm{I}}_{\phi}}d_q\widetilde{\mathcal{I}}_{-{\theta}}(x)\\
&=& e^{qt}d_q\Lambda_{{\rm{I}}_{\phi}}\widetilde{\mathcal{I}}_{-{\theta}}(x)=e^{qt}d_q\widetilde{\mathcal{I}}_{\Psi}(x),
 \end{eqnarray*}
 which shows, since it is positive,  that $d_q\widetilde{\mathcal{I}}_{\Psi}, q>0,$ is a $q$-invariant function for the  semigroup $P^{\Psi}$. Then, noting that $\widetilde{\mathcal{I}}_{\Psi}(0)=1$, an  application of  Dynkin's formula to the bounded stopping ${\rm{T}}_{\Psi}(t)=\inf({\rm{T}}_{\Psi},t)$, yields, for all $t\geq0$, that
 \begin{eqnarray*}%\label{eq:LTT}
e^{-qt}P^{\Psi}_{{\rm{T}}_{\Psi}(t)} d_q\widetilde{\mathcal{I}}_{\Psi}(0)&=&1
 \end{eqnarray*}
where $(P^{\Psi}_{{\rm{T}}_{\Psi}(t)})_{t \ge 0}$ is the semigroup corresponding
to the process $(\rm{X}_{{\rm{T}}_{\Psi}(t)})_{t \ge 0}$.
  Differentiating term by term the series $\widetilde{\mathcal{I}}_{\Psi}$,  we observe that $x\mapsto d_q\widetilde{\mathcal{I}}_{\Psi}(x)$ is, for all $q>0$, non-decreasing on $\R^+$. This allows us to invoke a dominated convergence argument, while combined with the absence of positive jumps of $\mathrm{X}(\Psi)$, which  entails that $\rm{X}_{{\rm{T}}_{\Psi}}(\Psi)=1$ almost surely,  gives the first identity in \eqref{eq:LTT}. The rest of the statement, that is,  ${\rm{T}}_{\Psi}$  is a positive self-decomposable variable is justified  in \cite{Patie}.
\end{proof}
 We are ready to state the following result that identifies the first random variable of our pair of solutions to the van Dantzig problem.
 \begin{lemma}\label{lem:BT}
Let  $\Psi \in \N$ with $0\leq \theta<1$ and ${\rm{B}}=({\rm{B}}_t)_{t\geq 0}$ be a standard Brownian motion independent of $\rm{X}(\Psi)$. Then, the random variable  $\overline{{\rm{D}}}_{\Psi}=\sqrt{2}{\rm{B}}_{{\rm{T}}_{\Psi}}$ is  real-valued, symmetric, and  infinitely divisible. Moreover,   for any $t\in \R$,
 \begin{equation}
 \mathcal F_{\overline{{\rm{D}}}_{\Psi}}(t)=\frac{1}{\mathcal{I}_{\Psi}(t)}.
 \end{equation}
   \end{lemma}
   \begin{proof}
  According to Theorem \ref{thm:patie}, ${\rm{T}}_{\Psi}$ is self-decomposable and hence infinitely divisible.
 Since ${\rm{B}}$ is a
symmetric L\'evy process, we have, by Bochner subordination, see \cite[Theorem 30.1]{Sato}, that $\overline{{\rm{D}}}_{\Psi}$ is  real-valued, symmetric, and  infinitely divisible. Next,  using the independence of ${\rm{B}}$ and $\rm{X}(\Psi)$ and hence of ${\rm{T}}_{\Psi}$, one gets, for any $t\in \R$,
 \begin{eqnarray*}
\mathcal F_{\overline{{\rm{D}}}_{\Psi}}(t)&=&\mathcal F_{{{\rm{T}}_{\Psi}}}\left(it^2\right)=\frac{1}{\mathcal{I}_{\Psi}(t)}
 \end{eqnarray*}
 where the last line follows from \eqref{eq:LTT}.
   \end{proof}
To characterize the second random variable, we need to introduce  the random variable ${\rm{J}}_{\nu}$, $\nu \in (-\infty,\frac12]$, whose law is, {{ when $\nu <-\frac12$ }}, absolutely continuous with a density $f_{{\rm{J}}_{\nu}}$ given by
\begin{equation}\label{eq:defJ}
f_{{\rm{J}}_{\nu}}(x)=\frac{2^{2\nu} \Gamma(1-{\nu})}{\Gamma(\frac{1}{2})\Gamma(\frac{1}{2}-{\nu})}\,(4-x^2)^{{-\nu}-\frac{1}{2}}\mathbb{I}_{\{|x|\leq 2\}}
\end{equation}
 {{and, otherwise for $\nu =-\frac12$, ${\rm{J}}_{\nu}$  has the Bernoulli distribution $\frac{1}{2}(\delta_2+\delta_{-2})$}}.
We note that the characteristic function corresponding to the random variable
${\mathrm J}_\nu$ is the function
$\mathcal J_{-\nu}(t) = \Gamma(1-\nu) t^\nu J_{-\nu}(2t)$,
see \cite[p.\ 38]{Bateman}.
We are now ready to state the following which defines  a new class of random variables indexed by the whole set of Bernstein  functions, and, when restricted to the subet $\Bd$ gives the other set of van Dantzig variables.
 \begin{lemma}\label{lem:IJ}
Let $\phi  \in \Be $ and
 define, for any {{$\nu \leq-\frac12$}},
 the random variable
\[
{\rm{D}}_{\phi,\nu}=\sqrt{{\rm{I}}_{\phi}} \times {\rm{J}}_{\nu}
\] where  ${\rm{J}_{\nu}}$ is chosen independent of  ${\rm{I}}_{\phi}$. Then,
${\rm{D}}_{\phi,\nu}$ is a symmetric random variable  taking values  in the possibly infinite interval   $\left(-\frac{2\sqrt{2}}{\sigma},\frac{2\sqrt{2}}{\sigma}\right)$ (we use the convention $\frac{1}{0}=\infty$). Moreover, all its even  moments exist and are given, {{for any $\nu < -\frac12$ (resp.~$\nu = -\frac12$)}}, by
 \begin{equation}\label{eq:moD}
 \mathcal{M}_{{\rm{D}}_{\phi,\nu}}(2n)=\frac{\Gamma(2n+1)\Gamma(1-{\nu})}{W_{\phi}(n+1)\Gamma(n+1-{\nu})}  \quad (\textrm{resp.~} =2^{2n}), \quad n= 0,1,\ldots,
 \end{equation}
 and, for any $t\in \R$, we have
 \begin{equation}\label{eq:FTD}
 \mathcal F_{{\rm{D}}_{\phi,\nu}}(t)=\mathcal{J}_{\phi,\nu}(t)
 \end{equation}
 where the entire function $\mathcal{J}_{\phi,\nu}$ is defined as in \eqref{eq:defIpz}.
 Finally,  the law of ${\rm{D}}_{\phi,\nu}$ is absolutely continuous with a density $f_{{\rm{D}}_{\phi,\nu}}$ which is continuous on $\R$ and $f_{{\rm{D}}_{\phi,\nu}} \in  C_0^{\infty}(\R\setminus\{0\})$ (resp.~$f_{{\rm{D}}_{\phi,\nu}} \in C_0^{p}(\R\setminus\{0\})$, where $p=\lceil\frac{2}{\sigma^2}\left(\overline{\mu}(0^+)+\phi(0)\right)-\nu-\frac{1}{2}\rceil$) if $\sigma^2=0$ or $\overline{\mu}(0^+)=\infty$ (resp.~otherwise), such that   $f_{{\rm{D}}_{\phi,\nu}}(x)=f_{{\rm{D}}_{\phi,\nu}}(-x)$, and, for any $n\in \N$ (resp.~$n=0,\ldots, p$),  $x>0$ and $a>\frac{1}{2}+n$,
\begin{equation}\label{eq:FTD}
 f^{(n)}_{{\rm{D}}_{\phi,\nu}}(x)=(-1)^n\frac{\Gamma(1-{\nu})}{2\pi i}\int_{a-i\infty}^{a+i\infty}x^{1-2z} \frac{\Gamma(z)}{\Gamma(z-n)} \frac{\Gamma(2z-1-2n)}{W_{\phi}(z-n)\Gamma(z-{\nu}-n)}dz.
 \end{equation}

%\textcolor{blue}{this is incomplete.. mellin transform and Mellin-Barnes representation of the density, asymptotic}

 %where $(\eta_k)_{k\geq 1}$ is a sequence of  iid Laplace random variables, and, $(z_k)_{k\geq 1}$ are the positive zeros of $J_{0}(z)$ the Bessel function of the
%first kind of order $0$.
   \end{lemma}
 %  \begin{remark}
%Note that since $\Phi$
%{ INCOMPLETE SENTENCE HERE!}
%   \end{remark}
    \begin{proof}
First, note that the symmetry property of ${\rm{J}}_{\nu}$ entails the one of ${\rm{D}}_{\phi,\nu}$. Since the support of ${\rm{J}}_{\nu}$ is $[-2,2]$ (or {{$\{-2,2\}$ when $\nu=\frac12$}}) and the one of ${\rm{I}}_{\phi}$ is $\left[0,\frac{2}{\sigma^2}\right]$, see e.g.~\cite[Theorem 2.4]{PS-B}, we deduce readily the one of  ${\rm{D}}_{\phi,\nu}$. Being symmetric, only the even moments are non-zero and are given, for any  $n= 0,1,\ldots$, by
\begin{equation*}
\mathcal{M}_{{\rm{D}}_{\phi,\nu}}(2n)=\mathcal{M}_{{\rm{J}}_{\nu}}(2n)\mathcal{M}_{{\rm{I}}_{\phi}}(n)
\end{equation*}
where we used that the random variables are independent. Using the identity \eqref{eq:moX}
and, observing, from the duplication formula of the gamma function, that, for {{$\nu <-\frac12$}},
\begin{eqnarray}\label{eq:MD}
% \nonumber % Remove numbering (before each equation)
\int_{-2}^{2}x^{2n}f_{{\rm{J}}_{\nu}}(x)dx
   &=& \frac{\Gamma(1-{\nu})}{\Gamma(\frac{1}{2})\Gamma(\frac{1}{2}-{\nu})}2^{2n}\, \int_{0}^{1}y^{n-\frac12}(1-y)^{{-\nu}-\frac{1}{2}}dy \nonumber \\
  &=& 2^{2n}\frac{\Gamma(1-{\nu})\Gamma(n+\frac{1}{2})}{\Gamma(\frac{1}{2})\Gamma(n+1-{\nu})}\nonumber \\
  &=& \frac{\Gamma(1-{\nu})(2n)!}{\Gamma(n+1-{\nu})n!}
\end{eqnarray}
we derive easily the identity \eqref{eq:moD}. The case $\nu=-\frac12$ follows easily and from the identity before \eqref{eq:MD}, the expression of moment is continuous (from below) in $\nu$.
Next, since, see e.g.~\cite[Theorem 1]{Jurek},  for any $t\in\R$,
\begin{equation}
 \mathcal F_{{\rm{J}}_{\nu}}(t)=\mathcal{J}_{-\nu}(t)=\sum^{\infty}_{n=0} \frac{(-1)^{n}\Gamma(1-\nu)}{\Gamma(n+1-\nu)n!}t^{2n}
 \end{equation}
and recalling, from \cite{PS-B}, that the law of ${{\rm{I}}_{\phi}}$ is absolutely continuous,  the independence of ${{\rm{I}}_{\phi}}$ and ${\rm{J}}_{\nu}$ yields, for any $t\in\R$, that
\begin{eqnarray*}
  \mathcal F_{{\rm{D}}_{\phi,\nu}}(t)&=& \int_0^{\infty} \mathcal F_{{\rm{J}}_{\nu}}(\sqrt{x}t)f_{{\rm{I}}_{\phi}}(x)dx \\ &=& \int_0^{\infty}\mathcal{J}_{\nu}(\sqrt{x}t)f_{{\rm{I}}_{\phi}}(x)dx
  \\&=&\sum^{\infty}_{n=0} \frac{(-1)^{n}\Gamma(1-\nu)\mathcal{M}_{{\rm{I}}_{\phi}}(n) }{\Gamma(n+1-\nu)n!}t^{2n}  \\&=& \mathcal{J}_{\phi,\nu}(t)
 \end{eqnarray*}
where to justify the interchange of sums in the third equality, we resort to a classical Fubini's argument whose details are described in \cite[Section 1.77]{Tit}, and relies on the fact that the series $ \mathcal{J}_{\phi,\nu}$ is absolutely convergent  on $\R$.
For the last identity we used \eqref{eq:moX}.
Next, it is not difficult to see that the computation \eqref{eq:MD} extends to any complex $z$ such that $\Re(z)>-\frac{1}{2}$, and, from \cite[Theorem 2.4]{PS-B}, the expression \eqref{eq:moX} also extends (at least) to $\Re(z)>-1$, to get, for $\Re(z)>-\frac{1}{2}$,
\begin{equation*}
  \mathcal{M}_{{\rm{D}}^2_{\phi,\nu}}(z)=\frac{\Gamma(2z+1)\Gamma(1-{\nu})}{W_{\phi}(z+1)\Gamma(z+1-{\nu})},
\end{equation*}
which proves \eqref{eq:moD}.
Then, recalling the Stirling formula of the gamma function, for any $a>0$ and $|b|$ large,
\begin{equation}
	\left|\Gamma(a+ib)\right| \sim C e^{-|b|\frac{\pi}{2}}|b|^{a-\frac{1}{2}}
	\end{equation}
where,
 here and below,
$C>0$ is a generic constant, we deduce, that, for any  $a>-\frac{1}{2}$ and $|b|$ large,
\begin{equation*}
	\left|\frac{\Gamma(2a+1+2ib)}{\Gamma(a+1+ib)\Gamma(a+1-\nu+ib)}\right| \sim C |b|^{\nu-\frac{1}{2}}.
	\end{equation*}
This combines with \eqref{eq:condClass11} when $\sigma^2>0$ and $\overline{\mu}(0^+)<\infty$, gives for any $\epsilon>0$ and $a>-\frac{1}{2}$,
\begin{equation}
	\left|\mathcal{M}_{{\rm{D}}^2_{\phi,\nu}}(a+ib)\right| \leq C |b|^{-\frac{2}{\sigma^2}\left(\overline{\mu}(0^+)+\phi(0)\right)+\nu-\frac{1}{2}+\epsilon}.
	\end{equation}
We deduce that the mapping 	$b\mapsto |b|^p\left|\mathcal{M}_{{\rm{D}}^2_{\phi,\nu}}(a+ib)\right| $ is integrable on $\R$ whenever $p<\frac{2}{\sigma^2}\left(\overline{\mu}(0^+)+\phi(0)\right)-\nu-\frac{1}{2}$. Then, invoking classical results on Mellin inversion, see e.g.~\cite[Section 1.7.4]{Patie-Savov}, we get  the Mellin-Barnes representation of $f_{{\rm{D}}^2_{\phi,\nu}}$ which takes the form for any $x>0$ and $a>\frac12$,
\begin{equation*}
 f_{{\rm{D}}^2_{\phi,\nu}}(x)=\frac{\Gamma(1-{\nu})}{2\pi i}\int_{a-i\infty}^{a+i\infty}x^{-z} \frac{\Gamma(2z-1)}{W_{\phi}(z)\Gamma(z-{\nu})}dz.
 \end{equation*}
 Next, with $\nu < -\frac12$,  since $f_{{\rm{D}}_{\phi,\nu}}$ is symmetric, we have that $f_{{\rm{D}}_{\phi,\nu}}(x)=xf_{{\rm{D}}^2_{\phi,\nu}}(x^2), x>0$, from where we deduce the Mellin-Barnes  representation of $f_{{\rm{D}}_{\phi,\nu}}$. Next,   since the mapping $z\mapsto \frac{\Gamma(2z-1)}{W_{\phi}(z)\Gamma(z-{\nu})}$ is meromorphic on $(a',a'+\frac{1}{2}), 0<a'<\frac12$ with a simple pole at $\frac{1}{2}$, an application of the residues theorem yields
\begin{equation*}%\label{eq:FTD}
 f_{{\rm{D}}_{\phi,\nu}}(x)=\frac{\Gamma(1-{\nu})}{2W_{\phi}(\frac12)\Gamma(\frac12-{\nu})}+\frac{\Gamma(1-{\nu})}{2\pi i}\int_{a'-i\infty}^{a'+i\infty}x^{1-2z} \frac{\Gamma(2z-1)}{W_{\phi}(z)\Gamma(z-{\nu})}dz
 \end{equation*}
 from where we easily conclude, as $a'<\frac12$, that $\lim_{x\downarrow 0 }f_{{\rm{D}}_{\phi,\nu}}(x)=\frac{\Gamma(1-{\nu})}{2W_{\phi}(\frac12)\Gamma(\frac12-{\nu})}$. Moreover, by Mellin inversion, see again \cite[Section 1.7.4]{Patie-Savov}, and by symmetry, we get that $f\in C_0(\R)$. From the same reference and by a similar reasoning, we obtain the expression and the smoothness properties for the successive derivatives for all $n\leq p$.
The cases $\sigma=0$ or $\overline{\mu}(0^+)=\infty$ follow easily by means of a similar reasoning and using  \eqref{eq:condClass1} in place of \eqref{eq:condClass11}.
   \end{proof}

   We continue our program with the following observation  which is the key step in proving Proposition \ref{prop:LukMap} for the $L^{(2)}$ mapping. To state it, we recall that $X(\gamma)$, $\gamma\in \R$, is the $\gamma$-length-biased random variable of a non-negative random variable $X$, if its $\gamma^{th}$ moment $m_{\rm X}(\gamma)$ is finite and
   \begin{equation*}
F_{{\rm X}(\gamma)}(dx) = \frac{x^{\gamma}}{m_{\rm X}(\gamma)} F_{\rm X}(dx).
\end{equation*}
This notion was used in the proof of the Lukacs mapping  defined in \eqref{it:dif}. %above, since one easily checks that $L^{(2)} \mathcal{F}_{X} (t)= \mathcal{F}_{X(2)}(t)$.
We are now ready to state the following. %To prove that the mapping $L^{(2)} \mathcal{J}_{\Psi} \in \De$, we need the following interesting observation.
   \begin{proposition}
 For any $\phi  \in \Be $, we have $\mathcal{T}_1\phi(u)=\frac{u}{u+1}\phi(u+1) \in \Be$. Moreover, with the notation of \eqref{eq:moX}, we have, for any $t\in \R$,
  \[\mathcal{F}_{{\rm{I}}_{\mathcal{T}_1\phi}} (t)= \mathcal{F}_{{\rm{I}}_{\phi}(1)}(t)\]
  and, for any twice continuously differentiable function $\mathcal{F}_X \in \Pp_+$,
   \begin{equation}\label{eq:comMap}
     L^{(2)}  \Lambda_{{\rm{I}}_{\phi}} \mathcal{F}_X= \Lambda_{{\rm{I}}_{\phi}(1)} L^{(2)}\mathcal{F}_X
   \end{equation}
where $\Lambda_{{\rm{I}}_{\phi}} f(t) = \E\left[f(t\sqrt{{\rm{I}}_{\phi}})\right]$ and ${\rm{I}}_{\phi}$ is chosen independent of $X$. In particular, for any $k=1,2,\ldots,$ and,  any $t\in \R$, with $\nu<\frac12$,
\begin{equation} \label{eq:L2F}
L^{(2)}_k \mathcal{F}_{{\rm{D}}_{\phi,\nu}}(t)=  \mathcal{F}_{{\rm{J}}_{\nu}(2k)\sqrt{{\rm{I}}_{\phi}(k)}}(t)=\mathcal{F}_{{\rm{D}}_{\phi,\nu}(2k)}(t).
\end{equation}
   \end{proposition}
\begin{proof}
The first statement can be found in \cite{CKP}. Then, observing that $W_{\mathcal{T}_1\phi}(n+1)= \prod_{k=1}^{n}\frac{k}{k+1}\phi(k+1)=\frac{W_{\phi}(n+2)}{W_{\phi}(1)(n+1)}$, we deduce, from \eqref{eq:moX}, that, for any  integer $n$,
\begin{equation}
\mathcal{M}_{{\rm{I}}_{\mathcal{T}_1\phi}}(n)=\frac{n!}{W_{\mathcal{T}_1\phi}(n+1)}=\frac{(n+1)!W_{\phi}(1)}{W_{\phi}(n+2)}=\frac{\mathcal{M}_{{\rm{I}}_{\phi}}(n+1)}{\mathcal{M}_{{\rm{I}}_{\phi}}(1)} =\mathcal{M}_{{\rm{I}}_\phi(1)}(n) \end{equation}
which yields the second claim since these random variables are moment determinate, see \cite{CPY}, and the characteristic function uniquely determines the law of a random variable.
Next, as $\mathcal{F}_X$ is twice continuously differentiable, one gets, for any $t\in \R$,
   \begin{eqnarray*}%\label{eq:comMap}
     L^{(2)}\Lambda_{{\rm{I}}_{\phi}} \mathcal{F}_X(t)&=&L^{(2)} \mathcal{F}_{X\sqrt{{\rm{I}}_{\phi}}}(t) \\&=&\frac{\mathcal{F}^{(2)}_{X\sqrt{{\rm{I}}_{\phi}}}(t)}{\mathcal{F}^{(2)}_{X\sqrt{{\rm{I}}_{\phi}}}(0)}  \\&=&\frac{1}{\mathcal{F}^{(2)}_{\sqrt{{\rm{I}}_{\phi}}}(0)}  \int_{0}^{\infty} \frac{\mathcal{F}^{(2)}_X(t\sqrt{x})}{\mathcal{F}^{(2)}_{X}(0)}xf_{{\rm{I}}_{\phi}}(x)dx\\&=&  \Lambda_{{\rm{I}}_\phi(1)} L^{(2)}\mathcal{F}_X(t)
   \end{eqnarray*}
 where for the third equality, we used that $X$ and ${\rm{I}}_{\phi}$ are independent and for the last one that $\mathcal{F}^{(2)}_{\sqrt{{\rm{I}}_{\phi}}}(0)=\int_{0}^{\infty} xf_{{\rm{I}}_{\phi}}(x)dx=\frac{1}{\phi(1)}$.   To prove the last relation, we recall that ${\rm{D}}_{\phi,\nu}=\sqrt{{\rm{I}}_{\phi}} {\rm{J}}_{\nu}$, where  ${\rm{J}_{\theta}}$ is chosen independent of  ${\rm{I}}_{\phi}$ and thus resorting to the commutation type relation \eqref{eq:comMap}, one gets, for any $t\in \R$,
    \begin{eqnarray*}
  L^{(2)} \mathcal{F}_{{\rm{D}}_{\phi,\nu}}(t)= L^{(2)}  \Lambda_{{\rm{I}}_\phi} \mathcal{F}_{{\rm{J}}_{\nu}}(t)= \Lambda_{{\rm{I}}_{\phi}(1)} L^{(2)} \mathcal{F}_{{\rm{J}}_{\nu}}(t)=\Lambda_{{\rm{I}}_{\phi}(1)}\mathcal{F}_{{\rm{J}}_{\nu}(2)}(t)
  =\mathcal{F}_{{\rm{J}}_{\nu}(2)\sqrt{{\rm{I}}_{{\rm{I}}_{\phi}(1)}}}(t)
\end{eqnarray*}
which provides the claim for $k=1$. Then, an induction argument gives for any $k$,
    \begin{eqnarray*}
  L^{(2)}_{k+1}\mathcal{F}_{{\rm{D}}_{\phi,\nu}}(t)= L^{(2)} \circ L^{(2)}_{k} \mathcal{F}_{{\rm{D}}_{\phi,\nu}}(t)= L^{(2)} \mathcal{F}_{{\rm{J}}_{\nu}(2k)\sqrt{{\rm{I}}_{{\rm{I}}_{\phi}(k)}}}
  =\mathcal{F}_{{\rm{J}}_{\nu}(2k+2)\sqrt{{\rm{I}}_{{\rm{I}}_{\phi}(k+1)}}}(t)
\end{eqnarray*}
which  completes the proof.
\end{proof}

\subsubsection{End of the proof of Theorem \ref{thm:main}}
To complete the proof of Theorem \ref{thm:main}, we take $\Psi\in  \Nd$,  and, recall from Proposition \ref{prop:WH}, that there exists $\phi \in \Be_\D$ and $0\leq \theta\leq \frac{1}{2}$ such that $\Psi(u)=(u-\theta)\phi(u)$. Then combining lemmas \ref{lem:BT} and \ref{lem:IJ}, we obtain, using the notation of the latter lemma, that is ${\rm{D}}_{\phi,\theta}=\sqrt{{\rm{I}}_{\phi}} {\rm{J}}_{\theta}$,
\begin{equation*}
   \mathcal F_{\overline{{\rm{D}}}_{\Psi}}(it) \mathcal F_{{\rm{D}}_{\phi,\theta}}(t)=1, \quad t\in\R.
\end{equation*}
Hence $[\mathcal F_{\overline{{\rm{D}}}_{\Psi}},  \mathcal F_{{\rm{D}}_{\phi,\theta}}]$  form a van Dantzig pair, which completes the proof of the theorem after invoking Theorem \ref{thm:patie}.
\subsubsection{Proof of Theorem \ref{thm:pick}}
We simply sketch the proof of the first claim as the detailed arguments can be found in \cite{PKS}. Let $\phi\in\Be_{P_1}$. Then, observing that $W_{\Psi}(n+1)=n!W_{\phi}(n+1), n\geq0$, we write $F_{\phi}= \frac{1}{W_{\phi}}$, that is $F_{\phi}$ is a function which is  analytic and zero-free on the half-plane  $\Re(z)>z_1$, $z_1$ being the largest zero of $\phi$ which is simple, and, at least on this former half-plane, $F_{\phi}$ is solution to the recurrence equation \begin{equation}\label{eq:reF}
F_{\phi}(z+1)=\frac{1}{\phi(z)}F_{\phi}(z).
\end{equation} From this recurrence equation, we get that $z_1$ is also a simple zero of $F_{\phi}$, and $F_{\phi}$ admits an analytic extension to  $\Re(z)>\rho_1>z_2$ due to the interlacing property. Then, the $1$-separation property entails that $\rho_1=z_1-1$ is a simple pole of $\phi$, and, hence, $\rho_1$ is neither a zero nor a pole of  $F_{\phi}$, which, yields, by the recurrence \eqref{eq:reF}, that $F_{\phi}$ admits an analytic extension to $\Re(z)>z_2=\rho_2-1$. An induction argument gives  that $F_{\phi}$ is indeed an entire function with the sequence of simple zeros $(z_k)_{k\geq 1}$, which completes the proof of the first claim by invoking Laguerre theorem \cite[Theorem 4]{Runck} as $\mathcal{J}_{\Psi}(z)=\sum^{\infty}_{n=0} \frac{(-1)^nF_{\phi}(n+1)}{n!}z^{2n}$. For the last claim, we invoke, for instance, the Mittag-Leffler function presented in example \ref{ex:ML} which has non-real zeros.

\subsubsection{Proof of Proposition \ref{prop:LukMap}}
The case $p=1$ is simply the combination of Lemma \ref{lem:L1}, Proposition \ref{prop:WH} and Theorem \ref{thm:main}.
To prove the case $p=2$, we take $\Psi \in \Nd$ and recall that ${\rm{D}}_{\phi,\theta}=\sqrt{{\rm{I}}_{\phi}} {\rm{J}}_{\theta}$, where $\Psi(u)=(u-\theta)\phi(u)$.  Then,  using \eqref{eq:FTD}, Proposition \ref{prop:WH} and the identity  \eqref{eq:L2F}, one gets, for any $k=1,2,\ldots,$
    \begin{eqnarray}
  L^{(2)}_k\mathcal{J}_{\Psi}(t)= L^{(2)}_k \mathcal{F}_{{\rm{D}}_{\phi,\theta}}(t)= \mathcal{F}_{{\rm{J}}_{\theta}(2k)\sqrt{{\rm{I}}_{{\rm{I}}_{\phi}(k)}}}(t).
\end{eqnarray}
On the other hand, using \eqref{eq:propT2}, and recalling the identity $\mathcal{I}_{\Psi}(t)=\mathcal{J}_{\Psi}(e^{i\frac{\pi}{2}}t)$, one has
 \begin{eqnarray}
  L^{(2)}_k\mathcal{I}_{\Psi}=\mathcal{I}_{\overline{\mathcal{T}}^{2k}_{\frac{1}{2}}\Psi}.
\end{eqnarray}
However, since $\Psi \in \Nd$, then $\frac{1}{2}\geq \theta$ and, from Proposition \ref{prop:WH}\eqref{it:tb2}, the mapping  $u\mapsto \overline{\mathcal{T}}^{2k}_{\frac{1}{2}}\Psi(u)=\overline{\mathcal{T}}^{1}_{\frac{1}{2}}\circ \overline{\mathcal{T}}^{2k-1}_{\frac{1}{2}}\Psi(u) = \frac{u-\frac{1}{2}}{u+\frac{1}{2}}\overline{\mathcal{T}}^{2k-1}_{\frac{1}{2}}\Psi(u+\frac{1}{2}) \in \N$ with $ \overline{\mathcal{T}}^{2k}_{\frac{1}{2}}\Psi(\frac{1}{2})=0$. Then,  $\overline{\mathcal{T}}^{2k}_{\frac{1}{2}}\Psi$ fulfills the requirement of Theorem \ref{thm:patie} and hence
\begin{eqnarray*}%\label{eq:LTT}
\mathcal{F}_{{\rm{T}}_{{\overline{\mathcal{T}}^{2k}_{\frac{1}{2}}\Psi}}}(iu)&=&\frac{1}{\mathcal{I}_{\overline{\mathcal{T}}^{2k}_{\frac{1}{2}}\Psi}(\sqrt{u})}
 \end{eqnarray*}
 which concludes the proof after invoking Lemma \ref{lem:BT}.

\vspace*{1cm}


\begin{thebibliography}{10}

\bibitem{Ahlfors}
L.V.\ Ahlfors.
{\em Complex Analysis.}
McGraw Hill, New York, 1979.

\bibitem{Bart-Patie}
C.\ Bartholm\'e and P.\ Patie.
Tur\'an inequalities and complete monotonicity for a class of entire functions,
Anal. Math.  47,  no. 3, 507--527, 2021.

\bibitem{Bateman}
H.\ Bateman.
{\em Tables of Integral Transforms}, Vol.\ 1.
McGraw-Hill, New York, 1954.

\bibitem{Biane}
P.\ Biane. La fonction Z\^eta de Riemann et les Probabilit\'es.
In: {\em La Fonction Z\^eta}, ed.\ by N.\ Berline and Claude Sabahh.
Ed.\ \'Ecole Polytechn., 166--193, 2011.

\bibitem{BianeP}
P.\ Biane, J.\ Pitman and M.\ Yor.
Probabilistic interpretation of the Jacobi theta and the Riemann zeta functions, via Brownian excursions,
{\em Bull.\ AMS} {\bf 38}, 435--465, 2001.

\bibitem{Boch}
S. Bochner.
{\em Harmonic Analysis and the Theory of Probability},
University of California Press, Berkeley and Los Angeles, 1955.

\bibitem{Borel}
E.\ Borel.
{\em Le\c{c}ons sur les Fonctions Enti\`eres}, 2nd edition.
Gauthiers-Villars et fils, Paris, 1921.

\bibitem{Brei}
L.\ Breiman.
{\em Probability}, Addison-Wesley Publishing Co., London,  1968.

\bibitem{CPY}
P.\ Carmona, F.\ Petit and M.\ Yor.
On the distribution and asymptotic results for exponential functionals of L\'evy processes. Exponential functionals and principal values related to Brownian motion, 73--130, \newblock {\em Rev.\ Mat.\ Iberoamericana}, Madrid, 1997.

\bibitem{chafait}
D.\ Chafa\"{\i}. A probabilistic proof of the Schoenberg theorem.
Available at \href{http://djalil.chafai.net/blog/2013/02/09/a-probabilistic-proof-of-the-schoenberg-theorem/}{Libres pens\'{e}es d'un math\'{e}maticien ordinaire}, 2013.

\bibitem{CKP}
M.\ Chazal, A.E.\ Kyprianou and P. \ Patie.
On some new transformations of Laplace exponents of spectrally
negative  L\'evy processes,
\newblock  A lifetime of excursions through random walks and Lévy processes - a volume in honour of Ron Doney's 80th birthday,
 157--180, Progr. Probab., 78, Birkhäuser/Springer, Cham,  2021.

\bibitem{curtis}
J. H.\ Curtiss.
A Note on the Theory of Moment Generating Functions.
{\em Ann.\ Math.\ Statist.} {\bf 13}, no. 4, 430-433, 1942.

\bibitem{DB}
N.G.\ de Bruijn.
The roots of trigonometric integrals,
{\em Duke Math.\ J.} {\bf 17}, 197-226, 1950.

\bibitem{DimitrovBenCheikh09}
D.K.\ Dimitrov and Y. Cheikh.
\newblock Laguerre polynomials as {J}ensen polynomials of {L}aguerre-{P}\'olya
  entire functions.
\newblock {\em J.\ Comput.\ Appl.\ Math.} {\bf 233}, no.\ 3, 703--707, 2009.

\bibitem{zeroM}
M.M.\ Dzhrbashyan.
{\em  Integral Transforms and Representations of Functions in the
Complex Domain}, Nauka, Moscow, 1966  (in Russian).

\bibitem{Gasper}
G.\ Gasper and M.\ Rahman,
\newblock {\em  Basic Hypergeometric Series}, second ed.,
Encyclopedia of Mathematics and its Applications {\bf 96},
Cambridge University Press, Cambridge,
With a foreword by Richard Askey, 2004.

\bibitem{Goldo}
A.A.\ Gol'dberg and I.V.\ Ostrovs'kii.
The growth of entire ridge functions with real zeros,
\newblock {\em Math.\ Phys.\ Funct.\ Anal.}, no.\ V, 3--10, 156 Akad. Nauk. Ukrain.
SSR Fiz.-Tehn\. Inst.\ Nizkih Temperatur Kharkov, 1974 (in Russian).

%\bibitem
 %T. Hasebe and A. Kuznetsov. On free stable distributions. Electron. Commun. Probab. 19 (56) (2014) 1--12

\bibitem{Ono}
M. Griffin, K. Ono, L. Rolen, D. Zagier.
Jensen polynomials for the Riemann zeta function and other sequences,
\newblock {\em Proceedings of the National Academy of Sciences}, 116 (23):11103-11110, 2019.

\bibitem{Wid}
I.I.\ Hirschman and D.V.\ Widder,
{\em  The Convolution Transform},
Princeton University Press, Princeton, New Jersey, 1955.


\bibitem{Oeuvres}
C.\ Hermite, H.\ Poincar\'e and E.\ Rouch\'e.
{\em Oeuvres de Laguerre, Tome I.}
Gauthiers-Villars et fils, Paris, 1898.


\bibitem{Hinds}
W. E.\ Hinds.
Moments of complex random variables related to a certain class of  characteristic functions.
\newblock {\em J.\ Sankhy$\bar{a}$ Ser.\ A} {\bf 36}, no.\ 2, 219--222, 1974.



\bibitem{Ismail}
M.E.H.\ Ismail and D.H.\ Kelker.
Special functions, Stieltjes transforms and infinite divisibility.
\newblock {\em  SIAM J.\ Math.\ Anal.} {\bf 10}, no.\ 5, 884--901,  1979.

\bibitem{LYs}
G. Iurato. The early historical roots of Lee-Yang theorem.
  {\em arXiv: History and Philosophy of Physics}, 2014.



\bibitem{Jurek}
Z.J.\ Jurek. Generalized L\'evy stochastic areas and selfdecomposability.
\newblock {\em Statist.\ Probab\. Lett. } {\bf 64}, no.\ 2, 213--222, 2003.

\bibitem{Kallenberg}
O.\ Kallenberg.
\newblock {\em Foundations of Modern Probability}, Springer, 2002.

\bibitem{Hyp11}
S.I.\  Kalmykov and D.B.\ Karp.
Log-concavity and Tur\'an-type inequalities for the generalized hypergeometric function. {\em Anal.\ Math.} {\bf 43}, 567-580, 2017.

%{
%\bibitem{Khos}
%D.\ Khoshnevisan.
%Schoenberg's theorem via the law of large numbers.
%\href{https://arxiv.org/abs/math/0504603v2}{ArXiv:math/0504603v2 [math.PR]}
%}


\bibitem{Kilbas06}
A.A.\ Kilbas, H.M.\ Srivastava, and J.J.\ Trujillo.
\newblock {\em Theory and Applications of Fractional Differential Equations},
 volume 204 of {\em North-Holland Mathematics Studies}.
\newblock Elsevier Science B.V., Amsterdam, 2006.


\bibitem{Koz}
Y.\ Kozitsky.
Laguerre entire functions and the Lee-Yang property,
\newblock {\em Applied Mathematics and Computation} {\bf 141}, no.\ 1, 103--112, 2003.

%\bibitem{Kos}
%H. Kosaki: On infinite divisibility of positive definite functions arising from
%operator means.\newblock {\em  J. Funct. Anal. } 254, 84--108, 2008.

\bibitem{KF48}
H.L.\ Krall and O.\ Frink.
A new class of orthogonal polynomials: the Bessel polynomials.
{\em Trans. Amer. Math. Soc.} {\bf 65}(1): 100--115, 1948.

\bibitem{PKS}
A.\ Kuznetsov, P.\ Patie and M.\ Savov.
Zeros of some entire functions and hitting times of self-similar Markov processes,
working paper, 2020.

\bibitem{Kyprianou06}
A.E.\ Kyprianou.
\newblock {\em Introductory Lectures on Fluctuations of {L}\'evy Processes with
Applications}.
\newblock Universitext. Springer-Verlag, Berlin, 2006.








\bibitem{Mateu}
M.\ Kwa\'snicki.
A new class of bell-shaped functions.
\newblock {\em Trans.\ Amer.\ Math.\ Soc.} {\bf 373} no.\ 4, 2255--2280, 2020.



\bibitem{Lamperti}
J.\ Lamperti.
Semi-stable Markov processes.
\newblock {\em  Z.\ Wahrsch.\ verw.\ Geb.} {\bf 22}, 205-225, 1972.

\bibitem{LY52}
T.D.\ Lee and C.N.\ Yang.
Statistical theory of equations of state and phase transitions, II.
Lattice gas and Ising model.
{\em Physical Rev.\ (2)} {\bf 87}, 410-419, 1952.

\bibitem{Levin96}
B.~Y. Levin.
\newblock {\em Lectures on entire functions}, volume 150 of {\em Translations
  of Mathematical Monographs}.
\newblock American Mathematical Society, Providence, RI, 1996.
\newblock In collaboration with and with a preface by Yu. Lyubarskii, M. Sodin
  and V. Tkachenko, Translated from the Russian manuscript by Tkachenko.


\bibitem{Linnik77}
Ju.V.\  Linnik and I.V.\ Ostrovskii.
{\em Decomposition of Random Variables and Vectors}.
Amer.\ Math.\ Soc., Providence 1977.
(Translation from the Russian original  published by Nauka, Moscow, 1972.)

\bibitem{Loeffen}
R.\ Loeffen, P.\ Patie and M.\ Savov.
Extinction time of non-Markovian self-similar processes, persistence,
annihilation of jumps and the Fr\'echet distribution.
{\em  J.\ Stat.\ Phys.} {\bf 175}, no.\ 5, 1022-1041, 2019.

\bibitem{Lukacs}
E.\ Lukacs.
Contributions to a problem of D.\ van Dantzig.
\newblock {\em Teor.\ Veroyatnost.\ i Primenen.} {\bf 13} no.\ 1, 114-125, 1968.

\bibitem{Lukacs_book}
E.\ Lukacs.
\newblock {\em  Characteristic Functions}, second edition.
Hafner Publishing Co., New York, 1970.

\bibitem{New1}
C.M.\ Newman.
\newblock Fourier transforms with only real zeros.
\newblock {\em Proceedings of the American Mathematical Society}
{\bf 61} no.\ 2, 245--251, 1976.

\bibitem{New}
C.M.\ Newman and W.\ Wu.
Lee-Yang property and Gaussian multiplicative chaos.
\newblock {\em Commun.\ Math.\  Phys.} {\bf 369}, 153--170, 2019.

\bibitem{NewRH}
C.M.\ Newman and W.\ Wu.
 Constants of de Bruijn-Newman type in analytic number theory and statistical physics.
\newblock {\em Bull.\ Amer.\ Math.\ Soc.},  in press.
Published online April 19, 2019.


\bibitem{Pakes}
A.G. Pakes.
Lambert's W, infinite divisibility and Poisson mixtures,
\newblock {\em Journal of Mathematical Analysis and Applications,}
Volume 378, Issue 2,
Pages 480-492, 2011.

\bibitem{Patie}
P.\ Patie.
\newblock Infinite divisibility of solutions to some self-similar
integro-differential equations and exponential functionals of {L}\'evy processes.
\newblock {\em Ann.\ Inst.\ Henri Poincar\'e Probab.\ Stat.}
{\bf 45} no.\ 3, 667--684, 2009.

\bibitem{Patie1}
P.\ Patie.
Law of the absorption time of positive self-similar Markov processes, {\em Ann.
Probab.} 40(2): 765-787, 2012.


\bibitem{Patie-Sarkar}
P.\ Patie and R.\ Sarkar.
Weak similarity orbit of the log-Bessel semigroup on the Euclidean space.
\newblock {\em  submitted}, 2020.

\bibitem{Patie-Savov2}
P.\ Patie and M.\ Savov.
Extended factorizations of exponential functionals of L\'evy processes.
\newblock {\em  Electron.\ J.\ Probab.} {\bf 17}, paper no. 38, 22 pp, 2012.

\bibitem{PS-B}
P.\ Patie and M.\ Savov.
\newblock
Bernstein-gamma functions and exponential functionals of {L}\'{e}vy processes.
\newblock {\em Electron.\ J.\ Probab.} {\bf 23}, paper no.\ 75, 101, 2018.

\bibitem{Patie-Savov}
P.\ Patie and M.\ Savov.
Spectral expansion of non-self-adjoint generalized Laguerre semigroups,
\newblock {\em Mem.\ Amer.\ Math.\ Soc.},  272,  no. 1336, vii+182 pp., 2021.

\bibitem{PSZ}
P.\ Patie,  M.\ Savov and Y.\ Zhao.
\newblock Intertwining, excursion theory and Krein theory of strings for non-self-adjoint Markov semigroups,
\newblock {\em Ann.\ Probab.} {\bf 47}, no.\ 5, 3231--3277, 2019.

\bibitem{Pitman-Yor-80}
J.\ Pitman and M.\ Yor,
\emph{{{B}essel processes and infinitely divisible laws}},
Stochastic Integrals (Proc.\ Sympos.\ Univ.\ Durham, Durham, 1980),
Lecture Notes in Math.\ (D.\ Williams, ed.) {\bf 851}, Springer, Berlin,
285--370, 1981.

\bibitem{Polya27}
G.\ P\'olya.
\"Uber trigonometrische Integrale mit nur reellen Nullstellen.
{\em J.\ reine u.\ ang.\ Math.}  {\bf 158}, 6-18, 1927.

\bibitem{PolyaSchur}
G.\ P\'olya and J.\ Schur.
\"Uber zwei Arten von Faktorenfolgen in der Theorie der Algebraischen Gleichungen.
{\em J.\ Reine Angew.\ Math.} {\bf 144}, 89-113, 1914.
\href{http://resolver.sub.uni-goettingen.de/purl?PPN243919689_0144}{Digital version by Nieders\"achsische Staats- und Universit\"atsbibliothek G\"ottingen}.

\bibitem{Popov}
A.Yu.\ Popov, A.M.\ Sedletskii.
Zeros distribution of Mittag-Leffler functions.
Contemporary Mathematics. Fundamental Directions 40 (2011), 3--171 (in Russian); EN Transl. in: J. Math. Sci. 190, 209--409, 2013.

\bibitem{Rivero}
V. Rivero. Recurrent extensions of self-similar Markov processes and Cram\'er’s condition. \newblock {\em Bernoulli}, 11(3), 471--509, 2005.

\bibitem{RT}
 B. Rodgers and T. Tao (2020). The De Bruijn–Newman constant is non-negative. Forum of Mathematics, Pi, 8, E6. doi:10.1017/fmp.2020.6

\bibitem{Roy}
B.\ Roynette and M.\ Yor.
Couples de Wald infiniment divisibles.
Exemples li\'es \`a la fonction Gamma d’Euler et \`a la fonction Zeta de Riemann.  \newblock {\em Ann.\ Inst.\ Fourier} {\bf 55} no.\ 4, 1219-1283, 2005.


\bibitem{Runck}
 H.-J.\ Runckel.
Zeros of entire functions.
\newblock {\em Trans.\ Amer.\ Math.\ Soc.} {\bf 143}, 343--362, 1969.

\bibitem{Sato}
K.-I.\ Sato.
\newblock {\em L\'evy processes and infinitely divisible distributions},
  volume~68 of {\em Cambridge Studies in Advanced Mathematics}.
\newblock Cambridge University Press, Cambridge, 1999.
\newblock Translated from the 1990 Japanese original, Revised by the author.

%{
%\bibitem{Sch1938}
%I.J.\ Schoenberg.
%Metric spaces and completely monotone functions.
%{\em Ann.\ Math.} (second series) {\bf 39}, no.\ 4, 811--841, 1938.
%}

\bibitem{Schoenberg}
I.J.\ Schoenberg.
On Totally positive functions, Laplace integrals and entire functions of the
Laguerre-Polya-Schur type.
{\em Proc.\ Nat.\ Acad.\ Sci.} {\bf 33}, no.\ 1, 11--17, 1947.

\bibitem{SchillingSongVondracek10}
R.L.\ Schilling, R.\ Song, and Z.\ Vondra{\v{c}}ek.
\newblock {\em Bernstein functions}, volume~37 of {\em de Gruyter Studies in
  Mathematics}.
\newblock Walter de Gruyter \& Co., Berlin, 2010.
\newblock Theory and applications.

\bibitem{SG73}
B.\ Simon and R.B.\ Griffiths.
The $(\phi^4)_2$ field theory as a classical Ising model.
{\em Comm.\ Math.\ Phys.} {\bf 33}, 145--164, 1973.

\bibitem{Tao}
T. Tao et al., Upper bounding the de Bruijn-Newman constant, available at  \href{https://terrytao.wordpress.com/2018/01/24/polymath-proposal-upper-bounding-the-de-bruijn-newman-constant/}{Polymath}
(2018)

\bibitem{Tit}
E.C.\ Titchmarsh.
\newblock {\em The Theory of Functions}.
Oxford University Press, Oxford, 1939.

\bibitem{YL52}
C.N.\ Yang and T.D.\ Lee.
Statistical theory of equations of state and phase transitions, I.
Theory of condensation.
{\em Physical Rev.\ (2)} { 87}, 404--409, 1952.

\bibitem{WidL}
D.V. Widder. \newblock {\em The Laplace Transform}.   Princeton University Press, 1941.

\bibitem{Zhang}
R.~Zhang.
On complete monotonicity of certain special functions. (English summary)
\newblock {\em Proc. Amer. Math. Soc.} 146, no. 5, 2049–2062, 2018.
\end{thebibliography}
\end{document}